%% file: manifoldSharp.tex
\documentclass[3p, number,sort&compress, times,fleqn]{elsarticle}

\usepackage{amsmath}
\usepackage{amssymb}
\usepackage{mathtools}

\usepackage[colorlinks = true,
            linkcolor = blue,
            urlcolor  = blue,
            citecolor = blue,
            anchorcolor = blue]{hyperref}

\usepackage[utf8]{inputenc}
\usepackage{subfig}
\usepackage{graphicx}

\usepackage{footnote}

\usepackage{./csml}

\newcommand{\ADDEDTEXTMODETWO}[1]{\textcolor{csmlRed}{#1}}
\newcommand{\ADDEDTWO}[1]{
\ifmmode
  \text{\ADDEDTEXTMODETWO{$#1$}}
\else
  \ADDEDTEXTMODETWO{#1}
\fi}

\let\latexchi\chi
\makeatletter
\renewcommand\chi{\@ifnextchar_\sub@chi\latexchi}
\newcommand{\sub@chi}[2]{
  \@ifnextchar^{\subsup@chi{#2}}{\latexchi^{}_{#2}}%
}
\newcommand{\subsup@chi}[3]{
  \latexchi_{#1}^{#3}%
}
\makeatother

\begin{document}

\begin{frontmatter}

\title{Manifold-based isogeometric analysis basis functions with prescribed sharp features}

\author{Qiaoling Zhang}  
\author{Fehmi Cirak\corref{cor1}}  
\ead{f.cirak@eng.cam.ac.uk}

\cortext[cor1]{Corresponding author}

\address{Department of Engineering, University of Cambridge, Trumpington Street, Cambridge CB2 1PZ, UK}

\begin{abstract}
We introduce manifold-based basis functions for isogeometric analysis of surfaces with arbitrary smoothness, prescribed $C^0$ continuous creases and boundaries. The utility of the manifold-based surface construction techniques in isogeometric analysis was demonstrated in Majeed and Cirak (CMAME,  2017). The respective basis functions are derived by combining differential-geometric manifold techniques with conformal parametrisations and the partition of unity method.  The connectivity of a given unstructured quadrilateral control mesh in~$\mathbb R^3$ is used to define a set of overlapping charts. Each vertex with its attached elements is assigned a corresponding conformally parametrised planar chart domain in~$\mathbb R^2$ so that a quadrilateral element is present on four different charts.  On the collection of  unconnected chart domains, the partition of unity method is used for approximation. The transition functions required for navigating between the chart domains are composed out of conformal maps.  The necessary smooth partition of unity, or blending, functions for the charts are assembled from tensor-product B-spline pieces and require in contrast to earlier constructions no normalisation.  Creases are introduced across user tagged edges of the control mesh. Planar chart domains that include creased edges or are adjacent to the domain boundary require special local polynomial approximants in the partition of unity method. Three different types of chart domain geometries are necessary to consider boundaries and arbitrary number and arrangement of creases.  The new chart domain geometries are chosen so that it becomes trivial to establish local polynomial approximants that are always $C^0$ continuous across the tagged edges. The derived non-rational manifold-based basis functions  correspond to the vertices of the mesh and may have an arbitrary number of creases and prescribed smoothness.  This makes them particularly well suited for isogeometric analysis of Kirchhoff-Love thin shells with kinks, which require $C^1$ continuous basis functions that are $C^0$ continuous across the kinks.  We demonstrate the convergence and utility of the new basis functions with linear and nonlinear beam, plate and shell examples.
\end{abstract}

\begin{keyword}
isogeometric analysis  \sep manifolds \sep smooth basis functions  \sep partition of unity method \sep sharp features  \sep thin shells 
\end{keyword}

\end{frontmatter}

\newpage
\includecomment{editing}


\newpage

\input{introduction}

\input{review}

\input{creased}

\input{mechanics}

\input{examples}

\input{conclusions}

\input{appendix}


\bibliographystyle{elsarticle-num-names}
\bibliography{manifoldSharp.bib}

\end{document}

%% file: introduction.tex
%
\section{Introduction}  

Smooth approximation schemes for unstructured meshes are crucial for isogeometric design and analysis of parts with arbitrary topology. Until recently isogeometric analysis was dominated by NURBS basis functions, which is the prevailing technology in present CAD systems. To represent the bounding surface of parts with arbitrary topology CAD systems usually resort to trimmed NURBS. Trimming involves the computation of the intersection between spline surfaces with other surfaces or curves. The respective   nonlinear root-finding problems look deceptively simple but are extremely hard to robustly solve and lead to non-watertight geometries~\cite{patrikalakis2009shape, marussig2017review, xiaoXiao:2018}. In the analysis context, the non-watertight geometries obtained from trimming pose unique challenges. Without turning to trimming, unstructured meshes with extraordinary vertices, i.e. vertices with different than four attached patches inside the domain, are required to represent parts with arbitrary topology. In computer-aided design numerous techniques have been developed to deal with extraordinary vertices, including geometrically $G^k$ and parametrically $C^k$ continuous constructions and subdivision surfaces. The application and further development of these techniques is currently a very active area of research in isogeometric analysis, see e.g.~\cite{Cirak:2000aa, scott2014isogeometric, buchegger2016adaptively, sangalli2016unstructured, collin2016analysis, nguyen2016c, toshniwal2017multi, toshniwal2017smooth, kapl2017isogeometric,kapl2018construction, chan2018isogeometric, zhang2018subdivision}. Unfortunately, none of the techniques from computer-aided design seems to give optimal finite element convergence rates without further modifications, especially when applied to Kirchhoff-Love thin shells with arbitrary geometry. Hence, the search for easy to implement and optimally convergent schemes
is still open. In a complementary line of  research, there has been progress in isogeometric analysis of shells based on trimmed surfaces and weak enforcement of mechanical continuity conditions across patch boundaries, see e.g.~\cite{breitenberger2015analysis, guo2018variationally}.

Manifold techniques for mesh-based construction of $C^k$ continuous surfaces were first introduced in Grimm and Hughes~\cite{grimm1995modeling}. In contrast to most other smooth surface construction techniques, which essentially rely on glueing of surface patches along their edges, in manifold techniques a surface is created by blending of overlapping surface patches. The surface patches are first defined over unconnected planar chart domains in~$\mathbb R^2$ and subsequently mapped to~$\mathbb R^3$.  On the collection of the unconnected planar chart domains, the partition of unity method~\cite{melenk1996partition} is used to construct the smooth surface patches. To that end, on each planar chart domain a local polynomial approximant and a partition of unity, or blending, function is needed. In addition, transition functions are required in order to be able to navigate between the different chart domains.  In manifold techniques the planar chart domains, the transition functions, the local polynomial approximants and the blending functions can all be relatively freely chosen, which makes them extremely versatile. Based on the seminal work of Grimm and Hughes a number of complementary mesh-based manifold approaches have been proposed in geometric modelling~\cite{grimm1995modeling, navau2000modeling, ying2004simple, della2008construction, beccari2014rags}. Most of these approaches differ in the choice of the planar chart domains. For instance, in~\cite{grimm1995modeling} each vertex, edge and face of the mesh have a corresponding planar chart domain with a suitably chosen geometry. In~\cite{ying2004simple} and~\cite{della2008construction} only the vertices have a corresponding planar chart domain in the form of a conformally mapped star-shaped polygonal disk or a circular disk, respectively.  In~\cite{navau2000modeling} the chart domains are chosen similar to the characteristic map in subdivision surfaces~\cite{Peters:2008aa}. Each of the mentioned choices for the chart domains implies a corresponding choice for the transition functions.  
Following the  construction proposed in Ying and Zorin~\cite{ying2004simple}, Majeed and Cirak~\cite{majeedCirak:2016}  introduced a new set of manifold-based basis functions  which can yield high convergence rates in finite element analysis. From a finite element viewpoint, the manifold-based basis functions resemble spline basis functions in the sense that each basis function has a  local support and has one corresponding vertex. 

\begin{figure}
  \centering
    \subfloat[][Control mesh with tags] 
  {
  	\includegraphics[scale=1.2]{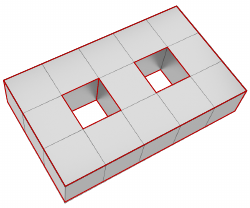}
  }
  \hfil
 \hspace{0.01\textwidth}
  \subfloat[][Manifold surface]
  {
  	\includegraphics[scale=1.2]{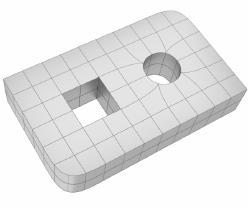}
  }
 \hspace{0.01\textwidth}
 \subfloat[][Deformed manifold surface]
 {
 	\includegraphics[scale=1.2]{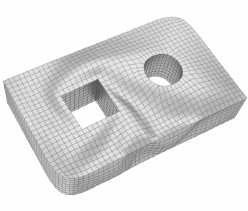}
 }  
  \caption{Isogeometric analysis of a genus 2 surface with creased $C^0$ continuous edges and sharp corners. (a) The control mesh with all the creased edges marked in red. (b)  The manifold surface with the faithfully reproduced creases. (c) The deformed manifold surface obtained with  (nonlinear) thin-shell finite element computation. In the computation the creases are modelled as rigid so that the angle across two adjacent surface pieces is maintained during the deformation. }
  \label{fig:blockIntro}
\end{figure}

In this paper, we derive new manifold-based basis functions for the isogeometric design and analysis of smooth surfaces with boundaries and with $C^0$ continuous sharp features, like creases and corners, by extending~\cite{majeedCirak:2016}.  The crease edges are tagged as such on the control mesh by the user and during the finite element analysis different mechanical continuity conditions can be imposed across the crease edges and along the boundary edges, see Figure~\ref{fig:blockIntro}. In mesh-based construction of surfaces using manifold techniques creases have previously been considered in Della Vecchia and J\"uttler~\cite{della2009piecewise} and boundaries in Tosun and Zorin~\cite{tosun2011manifold}.  The essential idea can be reduced to the choice of special local polynomials in the partition of unity approximation.  That is, the local polynomials on each chart domain have to consist out of several polynomial pieces that are~$C^0$ continuously connected across the crease edges. As in~\cite{majeedCirak:2016}, in our present construction each vertex and its attached elements have a corresponding star-shaped planar chart domain consisting of images of conformally mapped unit squares placed around a centre vertex. There is one chart domain per vertex and each quadrilateral element is present on four different chart domains. The choice of the special~$C^0$ continuous local polynomials can be simplified by slightly modifying the geometry of the chart domains.  The new chart domains are chosen such that they are rotationally symmetric with respect to the arrangement of crease edges. Or expressed differently, the crease edges must partition the chart domain into equiangular sectors. The~$C^0$ continuously connected polynomial pieces can subsequently be obtained by mapping a tensor-product basis, like the Lagrange or Bernstein basis, into each of the sectors of the chart domain.  The new chart domains are parameterised with a quasi-conformal map so that each element on it can have a different shape depending on the arrangement of the crease edges.  In contrast to the conformal map used in~\cite{majeedCirak:2016} the proposed quasi-conformal map is not angle-preserving but still provides  easily computable smooth transition functions. The proposed new construction leads to three different types of chart domains.  One of them is specifically designed to deal with crease arrangements which lead to creased sectors with concave corners.  As in~\cite{majeedCirak:2016} we assemble the smooth partition of unity blending functions from tensor-product B-spline segments defined over a unit square. However, similar to~\cite{tosun2011manifold}, the knot-interval and  the B-spline coefficients are chosen such that the assembled blending functions do not require normalisation. Hence, the new manifold-based basis functions derived in this paper are non-rational. 

The outline of this paper is as follows. In Section 2 the manifold-based basis functions introduced in~\cite{majeedCirak:2016} are briefly reviewed. 
The treatment of sharp features is discussed in Section 3. Depending on the arrangement of crease edges we distinguish between rotationally symmetric and asymmetric chart domains.  In addition, we consider the case of crease arrangements leading to concave sectors on a chart domain. Each of the three cases requires a different quasi-conformal map for parametrisation. Section 4 introduces the finite element analysis of thin shells with normal control along boundaries and across crease edges. In Section 5 the new manifold-based basis functions are applied to several Bernoulli beam and linear and nonlinear Kirchhoff-Love shell problems. The numerical convergence of $L_2$ and energy norm errors with decreasing element size is demonstrated.
%

%% file: review.tex
\section{Review of manifold-based basis functions}                          
%
In the following we briefly review the construction of univariate and bivariate manifold-based basis functions. The discussion is focused on their application in finite element analysis, so that the underlying manifold concepts from differential geometry and the partition of unity interpolation are only mentioned in passing. For a more comprehensive discussion on manifold-based basis functions and surface construction in geometry we refer to~\cite{majeedCirak:2016, ying2004simple, grimm1995modeling}.

\subsection{Univariate basis functions \label{sec:univariateBasis}}                          
\begin{figure}[t]
  \centering
	\includegraphics[scale=0.875]{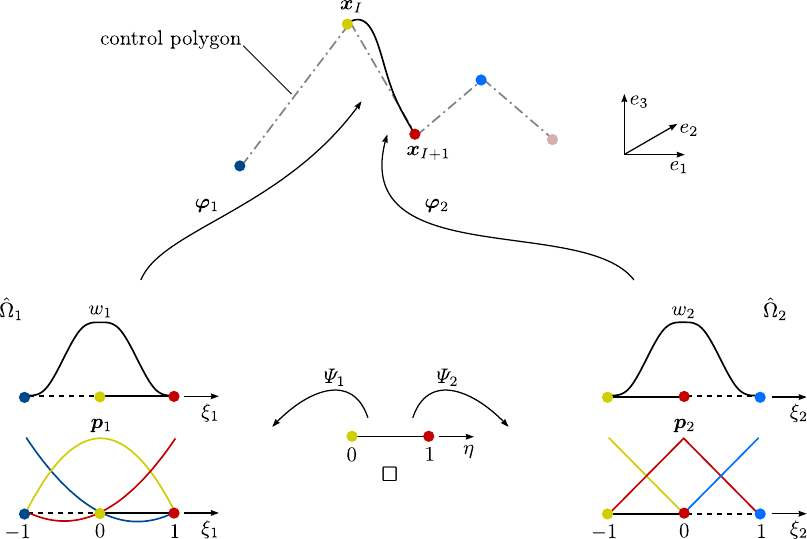}		
   \caption{Univariate manifold construction over a reference element~$\Box$ and the approximation of a given control mesh (dash-dotted). The reference element~$\Box$ maps into the two elements (solid lines) in the two chart domains~$\hat \Omega_1$ and~$\hat \Omega_2$. The local polynomial basis on~$\hat \Omega_1$ is a quadratic and on~$\hat \Omega_2$ it is a piecewise linear Lagrange basis. The blending functions~$w_1(\xi_1) = w_1(\varPsi_1(\eta))$ and~$w_2 (\xi_2) = w_2(\varPsi_2(\eta))$  sum up to one for~$\eta \in \Box$. The same vertex in different domains has the same colour. \label{fig:1dConstruction}}	
\end{figure}

%
It is instructive to first review the derivation of the univariate manifold-based basis functions.  We consider the dash-dotted control polygon in Figure~\ref{fig:1dConstruction} representing a part of a finite element mesh consisting of vertices \mbox{$\vec x_I \in \mathbb R^3$} and elements between consecutive vertices. Our aim is to derive the basis functions for a representative element~\mbox{$[ \vec x_I, \, \vec x_{I+1}] \in \mathbb R^3$,} as highlighted in Figure~\ref{fig:1dConstruction}. As in conventional finite elements we define a reference element~\mbox{$\Box \coloneqq [0, \, 1] \in \mathbb R$}  that will serve as an integration domain for evaluating the finite element integrals. In the manifold-based approach the  basis functions are obtained by smoothly blending local polynomials defined over several overlapping charts, or patches as they were called in~\cite{majeedCirak:2016}\footnote{We refrain in this paper from using the term patches because in computer-aided design literature patches denote what are  the elements in the finite element literature.}. In the following two chart domains~\mbox{$\hat \Omega_1 \coloneqq [-1, \, 1] \in \mathbb R$} and~\mbox{$\hat \Omega_2 \coloneqq  [-1, \, 1] \in \mathbb R$} are introduced for the element~$[ \vec x_I, \, \vec x_{I+1}]$. $\hat \Omega_1$ is associated with the vertex~$\vec x_I$ and its 1-neighbourhood and $\hat \Omega_2$ is associated with the vertex~$\vec x_{I+1}$ and its 1-neighbourhood. The 1-neighbourhood of a vertex is defined as the union of elements that contain the vertex. Hence, the  basis functions in the element~\mbox{$[ \vec x_I, \, \vec x_{I+1}]$} will be obtained by blending local polynomials defined over~$\hat \Omega_1$ and~$\hat \Omega_2$. 

As shown in Figure~\ref{fig:1dConstruction}, the reference element~$\Box$ maps to elements in both chart domains~$\hat \Omega_1$ and~$\hat \Omega_2$.  The respective maps are defined as  
\begin{align}
\begin{split}
		{\varPsi}_1  \colon & \eta \in \Box \mapsto  \xi_1 \in [0, \, 1]   \in \hat \Omega_1  \, , \\
		\varPsi_2  \colon  & \eta \in \Box \mapsto  \xi_ 2 \in [-1, \,  0] \in \hat \Omega_2 \, .
\end{split}
\end{align}
 In this specific univariate construction both maps are simple translations. The two maps imply a transition map between the two chart domains given by 
 \begin{equation}
 	  \varPsi_2 \circ \varPsi^{-1}_1 \colon \xi_1 \in \hat \Omega_1 \mapsto \xi_2 \in \hat \Omega_2 \, .
 \end{equation}
  We choose on each chart domain a polynomial approximant of the form
\begin{equation}
	f_j (\xi_j) = \vec p_j^\trans(\xi_j ) \vec \alpha_j  \quad  \text{with} \quad  j \in \{ 1, \, 2\} \, ,
\end{equation}
where the vectors~$\vec p_j(\xi_j)$ contain a local polynomial basis, like the power, Lagrange or Bernstein basis, and the vectors~$\vec \alpha_j$ contain their coefficients. In the present example a Lagrange basis is chosen, see~Figure~\ref{fig:1dConstruction}. The local basis  on the  chart~$\hat \Omega_1$ is quadratic and on the chart~$\hat \Omega_2$ it is piecewise linear. In addition to the polynomial approximants~$f_j(\xi_j)$, on each chart a smooth blending, or weight function~$w_j(\xi_j)$ is chosen. The two non-zero blending functions over the reference element~$\Box$ must satisfy the partition of unity property, i.e.,    
\begin{equation}\label{eq:partitionOfUnity}
	 \sum_{j=1}^2 w_j (\xi_j  )  = 1 \quad \forall \eta \in \Box  \quad \text{with}\quad \xi_1 = {\varPsi}_1(\eta) \, , \quad \xi_2 = {\varPsi}_2(\eta)  \, .
\end{equation}
\begin{figure}
  \centering
	\includegraphics[scale=0.875]{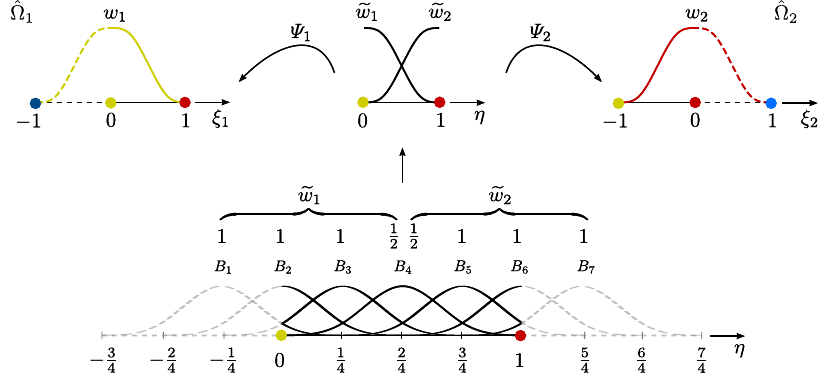}		
   \caption{Construction of smooth blending functions~$w_1(\xi_1)$ and~$w_2(\xi_2)$ as the linear combination of cubic B-splines defined over a parameter space with a knot-distance~$1/4$. The coefficients for~$\widetilde{w}_1(\eta)$ and~$\widetilde{w}_j(\eta)$ are indicated above the B-splines. \label{fig:1dWeightFun}}	
\end{figure}
The blending functions can be conveniently combined from suitably defined B-spline basis functions. In Figure~\ref{fig:1dWeightFun} the construction of $C^2$-continuous blending functions from cubic B-splines is illustrated. The uniform knot interval length of the B-splines is~$1/4$. The blending functions are defined as the linear combinations
\begin{align} \label{eq:weightFun}
\begin{split}
	\widetilde{w}_1(\eta) &= B_1(\eta) + B_2(\eta) + B_3(\eta) + \frac{B_4(\eta)}{2}  \, , \\
	\widetilde{w}_2(\eta) &= \frac{B_4(\eta)}{2} + B_5(\eta) + B_6(\eta) + B_7(\eta) \, ,
\end{split}	
\end{align}
and yield after mapping the blending functions 
\begin{equation}
 	w_j (\xi_j) = \widetilde{w}_j(\varPsi_j^{-1} (\xi_j))  \, .
\end{equation}
Evidently, both functions satisfy the partition-of-unity property~\eqref{eq:partitionOfUnity} given that a complete B-spline basis sums up to one. It is worth emphasising that the blending functions proposed here are polynomial in contrast to the rational ones in~\cite{ying2004simple, majeedCirak:2016}, see also~\ref{sec:weightsApp}. The tensor products of the new univariate blending functions yield the corresponding multivariate blending functions. Note that different from the chosen knot interval of 1/4, as proposed in~\cite{tosun2011manifold}, a knot interval  of 1/3 yields also blending functions which sum up to one without normalisation.

Finally, the approximant over the reference element~$\Box$ is obtained by blending the approximants over the two charts  
\begin{equation} \label{eq:approx}
	f(\eta) =   \sum_{j=1}^2 w_j (\xi_j  ) \vec p_j^\trans(\xi_j  ) \vec \alpha_j   \quad \text{with}\quad \xi_j = {\varPsi}_j(\eta) \, . 
\end{equation}
Due to the choice of the Lagrange basis for $\vec p_j(\xi_j)$ the coefficients $\vec \alpha_j$ can be interpreted as vertex coefficients. The coefficients $\vec \alpha_j$ of each chart are formally obtained with
\begin{equation} \label{eq:gather}
\vec \alpha_j = \vec P_j \vec f  \, , 
\end{equation}
where~$\vec f$ is the array of coefficients of all control polygon vertices and $\vec P_j$ is a gather matrix, filled with ones and zeros. This introduced in~\eqref{eq:approx} yields the manifold-based basis functions 
\begin{equation} \label{eq:approxN}
	f(\eta) =     \sum_{j=1}^2  \left  ( w_j (\xi_j  ) \vec p_j^\trans(\xi_j ) \vec P_j \right ) \vec f = \vec N^\trans (\eta ) \vec f   \quad \text{with}\quad \xi_j = {\varPsi}_j(\eta)\, ,
\end{equation}
where~$\vec{N}(\eta)$ is the array of non-zero basis functions over the considered element. The number of non-zero basis functions in~$\vec{N}(\eta)$  depends on the cardinality of the set of vertices in the two chart domains $\hat{\Omega}_1$ and $\hat{\Omega}_2$, and the specific local approximants on them. 

To investigate the smoothness of the basis functions~$\vec N(\eta)$ it is necessary to consider  their derivatives. According to definition~\eqref{eq:approxN}, their smoothness depends on the blending functions~$w_j(\xi_j)$, the local basis functions $\vec p_j(\xi_j)$ and the maps $\varPsi_j (\eta)$. To obtain $C^k$-continuous basis functions each of these functions has to be $k$-times differentiable over each chart domain; and,  in addition, the blending functions~$w_j(\xi_j)$  and their up to k-th derivatives must vanish at the chart domain boundaries, i.e.,
\begin{equation}
	 \frac{d^l w_j (-1) }{d \xi_j} = \frac{d^l w_j (1) }{d \xi_j} = 0 \quad \forall \,  l \leq k\, .
\end{equation}
Moreover, the differentiability of the mappings~$\varPsi_1(\eta)$ and~$\varPsi_2(\eta)$ requires that they satisfy at the centre of the chart domain 
\begin{equation}
	\frac{d^l \varPsi_1 (0) }{ d \eta^l} = \frac{ d^l  \varPsi_2 (1)  }{ d \eta^l}  \quad \forall \,  l \leq k\, .
\end{equation}

\begin{figure}[]
  \centering
   \subfloat[][Basis functions] 
  {
  	\includegraphics[scale=0.875]{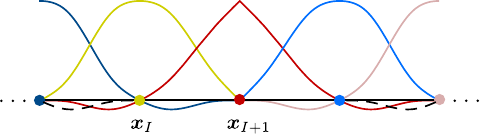}
  }
    \hspace{0.05\textwidth}
   \subfloat[][First derivatives] 
  {
  	\includegraphics[scale=0.875]{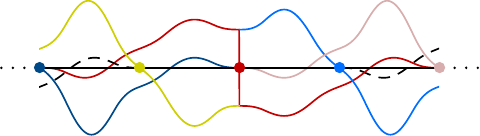}
  }
  \caption{Univariate basis functions and their derivatives over four elements. Five basis functions (solid) are non-zero over the two elements of the centre chart.}
  \label{fig:1dSfun}
\end{figure}

The smoothness of the basis functions~$\vec N(\eta)$ can be pointwise reduced by selecting a suitable local polynomial basis~$\vec p_j (\xi_j)$, as  in Figure~\ref{fig:1dConstruction} with a piecewise linear basis on chart~$\hat \Omega_2$. The resulting basis functions~$\vec N(\eta)$ and their first derivatives are plotted in Figure~\ref{fig:1dSfun}. The smoothness of the basis is reduced to $C^0$ at the vertex~$\vec x_{I+1}$. This has also an effect on the number of non-zero basis functions in an element. For instance, in the two elements~$[ \vec x_I\,, \vec x_{I+1} ]$  and~$ [ \vec x_{I+1}, \, \vec x_{I+2} ]$ with the linear piecewise $C^0$ continuous local basis functions~$\vec p_j(\xi_j)$ there are three non-zero basis functions. In elements with a quadratic local basis on each chart there are four non-zero basis functions. Finally, with the generated manifold-based basis functions points on the reference element~$\eta \in \Box$ are mapped onto the manifold according to
\begin{equation}
	\vec x (\eta) = \sum_{I} N_I(\eta) \vec x_I \,.
\end{equation}

\subsection{Bivariate basis functions}       
\label{sec:bivariateBasis}                   
\begin{figure}
	\centering 
	\includegraphics[scale=1.075]{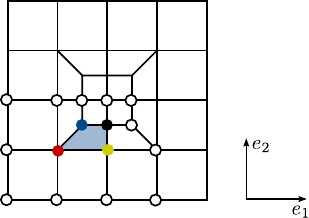}	
	\caption{An unstructured mesh with extraordinary vertices of valence $v \in \{3, \, 5 \} $. The shaded element is overlapped by four charts and the union of the four charts has 16 unique vertices in total. The reference element and the chart domains corresponding to the shaded element are shown in Figure~\ref{fig:manifoldPUM}. \label{fig:mesh}}
	\vspace*{5\floatsep}
%
	\centering 
	\includegraphics[scale=0.875]{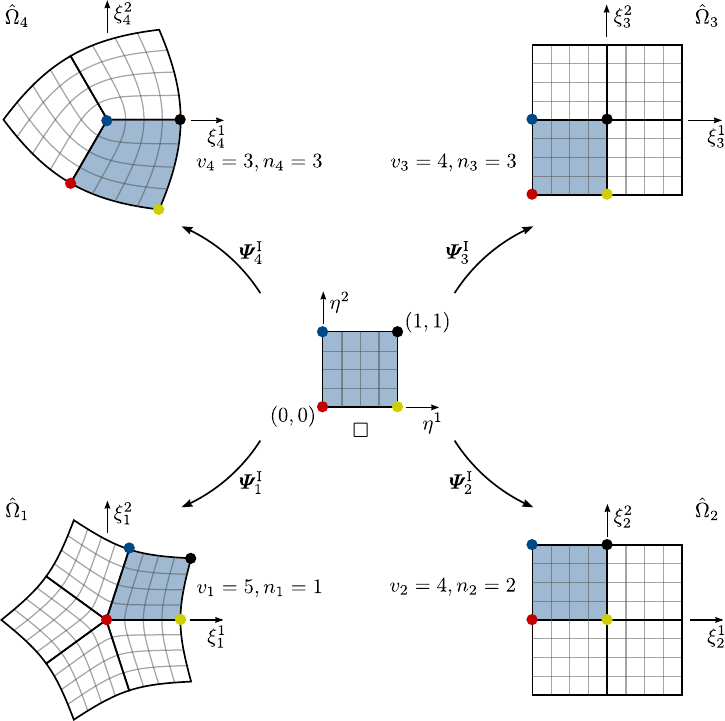}	
	\caption{The reference element $\Box$ (centre) and its four overlapping chart domains $\hat \Omega_j$ with $j\in \{1, \, 2, \, 3, \, 4\}$. The reference element~$\Box$ is mapped with $\vec \varPsi^\textrm{I}_j$ to the blue shaded faces in the four charts $\hat \Omega_j$. The variables~$n_j$ denote the face number of the shaded face and~$v_j$ the valence of the centre vertex. Notice that the parameter lines in the charts (with $\eta^1 = \text {const.}$ or $\eta^2= \text{const.}$)  are always orthogonal to the spoke edges which guarantees that the derivatives of the maps $ {\vec \varPsi}^\textrm{I}_j$  are  continuous over the entire chart domain. 	\label{fig:manifoldPUM}}
\end{figure}

%
The bivariate manifold-based basis functions provide smooth approximants even on unstructured surface meshes with extraordinary vertices. We consider the construction of the manifold-based basis functions  for a representative element in the quadrilateral finite element mesh shown in Figure~\ref{fig:mesh}. The reference element is now defined as the unit square \mbox{$\Box := [0, \, 1] \times [0, \, 1]  \in \mathbb R^2$}.  Similar to the univariate case each of the four vertices of the considered element and their 1-neighbourhoods have an associated chart domain \mbox{$\hat \Omega_j  \in \mathbb R^2$} with \mbox{$j \in \{ 1, \, 2, \, 3, \, 4 \}$}. The number of elements in a chart~$\hat \Omega_j$ depends on the number of elements~$v_j$ connected to the respective vertex, which is called the valence of the vertex. In quadrilateral meshes, the interior vertices with $v_j \neq 4$ are  the extraordinary  vertices. The basis functions will be obtained by smoothly blending local polynomials  defined over each of the four overlapping charts~$\hat \Omega_j$. 

As shown in~Figure~\ref{fig:manifoldPUM}, the reference element~$\Box$ maps to elements in the four chart domains~$\hat \Omega_j$ according to
\begin{equation}
	\vec \varPsi^\textrm{I}_j \colon \vec \eta = (\eta^1, \, \eta^2) \in \Box \mapsto \vec \xi = (\xi_j^1, \, \xi_j^2) \in \hat \Omega_j \, .
\end{equation}
The construction of these maps requires special care and will be detailed further below. The bivariate approximant  over the reference element $\Box$ is obtained by blending the approximants over the four charts 
\begin{equation} \label{eq:bivApprox}
	f(\vec{\eta}) = \sum_{j=1}^{4} w_j(\vec{\xi}_j )   f_j (\vec \xi_j)= \sum_{j=1}^{4} w_j(\vec{\xi}_j )  \vec p_j^\trans(\vec \xi_j  )  \vec \alpha_j  \quad \text{with}\quad \vec{\xi}_j = \vec{\varPsi}^\textrm{I}_j(\vec{\eta}) \, ,
\end{equation}	
where $w_j (\vec{\xi}_j )$ and $ \vec p_j(\vec \xi_j  )$ are the blending functions and the vector of local basis functions on the chart domain~$\hat \Omega_j$, respectively. As mentioned in Section~\ref{sec:univariateBasis}, the blending functions~$w_j (\vec{\xi}_j )$ are obtained by taking the tensor-product of univariate blending functions~\eqref{eq:weightFun}.  Next, the coefficients~$\vec \alpha_j$ are to be expressed in dependence of the vertex coefficients collected in the array~$\vec f$. The number of unique vertices in the union of the four charts~\mbox{$\cup_{j=1}^4 \hat \Omega_j$} is not fixed  and depends on the valences $v_j$ of the four vertices of the element. Hence, a least-squares projection is applied to express~\eqref{eq:bivApprox} in dependence of the vertex coefficients 
\begin{equation}
	f (\vec \eta) = \sum_{j=1}^4  \left (  w_j (\vec \xi_j) \vec{p}_j^\trans(\vec{\xi}_j) \vec{A}_j\vec{P}_j \right ) \vec{f}  \, ,
\end{equation}
where we used~\mbox{$\vec \alpha_j = \vec A_j \vec P_j \vec f$} on each chart~$\hat \Omega_j$. Here,~$\vec P_j$ is the gather matrix and~$\vec A_j$ the least-squares projection matrix. It is evident that the number of polynomial coefficients in~$\vec \alpha_j$  must not be greater than the number of vertices in the chart domain~\mbox{$\hat \Omega_j$.} It is worth mentioning that the projection matrix~$\vec A_j$ depends only on the valence~$v_j$ of the vertex and the chosen local basis function~$\vec p_j$ so that it can be precomputed and tabulated. Finally, the array of basis functions is defined with
\begin{equation} \label{eq:basisBiv}
	\vec{N}^\trans(\vec{\eta}) = \sum_{j=1}^{4} w_j(\vec{\xi}_j)\vec{p}_j^\trans (\vec{\xi}_j) \vec{A}_j\vec{P}_j  \, .
\end{equation}

As in the univariate case, the smoothness of the basis functions~$\vec N (\vec \eta)$ depends on the blending functions~$w_j(\vec \xi_j)$, the local basis functions~$\vec p_j(\vec{\xi}_j)$ and the mappings~$ \vec{\varPsi}^\textrm{I}_j(\vec{\eta}) $.  In the presented construction conformal maps are critical for the smooth parametrisation of  the chart domains~$\hat \Omega_j$ and the composition of the blending functions~$w_j(\vec \xi_j)$. Specifically, the map~$\vec{\varPsi}^\textrm{I}_j$ is composed of several conformal maps. 
To write conformal maps more succinctly, the coordinates $\vec \eta = (\eta^1, \, \eta^2 )$ of points in the reference element~$\Box$ are expressed as a complex number
\begin{equation}
	z = \eta^1 + i\eta^2 = |z|(\cos \phi + i\sin \phi) = |z|e^{i\phi} \, \quad \text{with} \quad |z| = \sqrt{(\eta^1)^2 + (\eta^2)^2} \quad \text{and} \quad \phi= \arctan(\eta^2/\eta^1)	 \, ,
\end{equation}
where $ |z|$ is the radius and  $\phi$ the angle in the complex plane.  In the complex plane the coordinates of the four corners of the reference element are given by 
\begin{equation}
	\begin{bmatrix} 
		z_1 & z_2 & z_3 & z_4
	\end{bmatrix}^\trans 
	= 
	\begin{bmatrix}
		0 + 0i &  1+0i &   1+1i & 0+1i
	\end{bmatrix}^\trans \, .
\end{equation}

\begin{figure}[]
	\centering 
	\includegraphics[scale=0.875]{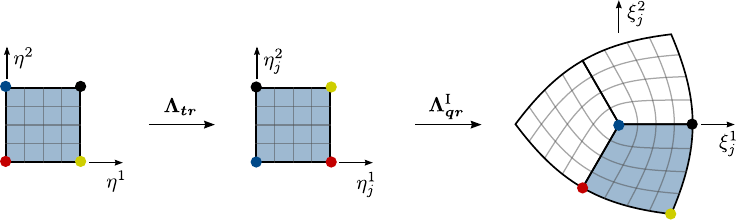}	
	\caption{The map from the reference element to a chart domain corresponding to a vertex with valence $v_j=3$ as denoted~$\vec \varPsi^{\textrm{I}}_4$  in Figure~\ref{fig:manifoldPUM}. The values $j=4$, $v_j = 3$, $n_j = 3$ are substituted into~\eqref{eq:linearMap} and~\eqref{eq:conformalMap} to obtain the expressions for $\vec {\Lambda_{tr}}$ and $\vec{\Lambda^\textrm{I}_{qr}}$.}
	\label{fig:manifoldMap}
\end{figure}

As illustrated in Figure~\ref{fig:manifoldMap}, the map~\mbox{$\vec{\varPsi}^\textrm{I}_j \colon \vec{\eta}  \in \Box \mapsto \vec{\xi}_j \in \hat \Omega_j$} is composed of two maps, i.e.,  
\begin{equation}\label{eq:composedMap}
	\vec{\varPsi}^\textrm{I}_j =  \vec{\Lambda_{qr}^\textrm{I}} \circ \vec{\Lambda_{tr}} \, .
\end{equation}
The auxiliary linear map $\vec{\Lambda_{tr}}$ is responsible for translating and rotating the reference element and is given by  
\begin{equation}
	\label{eq:linearMap}
	\vec{\Lambda_{tr}}(z;j) = (z - z_j) e^{-i \pi (j-1)/2} \, ,
\end{equation}
where $ (z - z_j)$ is a translation and $e^{-i \pi (j-1)/2}$ is a (rigid body) rotation. The quasi-conformal map~$\vec{\Lambda^\textrm{I}_{qr}}$ maps the reference element into a wedge-shaped domain and  is chosen as 
\begin{equation}
	\label{eq:conformalMap}
	\vec{\Lambda^\textrm{I}_{qr}}(z;v_j,n_j) = \frac{ |z|^{\beta}}{|z|^{4/v_j}}   z^{4/v_j}  e^{i 2\pi (n_j-1)/v_j} = |z|^\beta e^{i\phi 4/v_j}  e^{i 2\pi (n_j-1)/v_j} 	\, ,
\end{equation}
where $\beta$ is a free parameter, $v_j$ denotes the valence of the centre vertex and $n_j$ is the number of the sector onto which the reference element is mapped. According to the first expression in~\eqref{eq:conformalMap},  the map $\vec{\Lambda^\textrm{I}_{qr}}$ is composed of a   standard conformal map~$z^{4/v_j} $,  a scaling of the radius with~$|z|^{\beta-4/v_j}$ and a rotation~$ e^{i 2\pi (n_j-1)/v_j} $.   The parameter $\beta$ can be chosen relatively freely~\cite{ying2004simple}. For instance,  when $\beta = 4/v_j$ the scaling term drops and the map $\vec{\Lambda^\textrm{I}_{qr}}$ is composed of a conformal map and a rotation as in~\cite{majeedCirak:2016}. The second expression in~\eqref{eq:conformalMap} shows that the  radius $|z|$ will remain constant when $\beta=1$. Hence, in this paper we choose $\beta = 1$ to obtain a more uniform mapping $\vec{\varPsi}^\textrm{I}_j$ with smaller entries in its Jacobian matrix. As it will become clear, the choice~$\beta=1$ also simplifies the parametrisation of the asymmetric chart domains and chart domains with concave corners in Sections~\ref{sec:type2chart} and~\ref{sec:type3chart}.  See~\ref{sec:conformApp} for a review on conformal maps and a discussion on the choice of  $\beta$.

With the mapping $\vec{\varPsi}^\textrm{I}_j$ at hand, we can now evaluate the manifold-based basis functions defined in~\eqref{eq:basisBiv}. For instance, for a given integration point~$\vec \eta$ in the reference element, first its image~$\vec \xi_j = \vec \varPsi^\textrm{I}_j (\vec \eta)$ in each of the four overlapping charts is found and after that the blending functions~$w_j(\vec{\xi}_j)$ and local polynomials~$\vec{p}_j(\vec{\xi}_j)$ are evaluated.

%% file: creased.tex
\section{Creased bivariate manifold-based basis functions}                          
%
We introduce  domain boundaries and~$C^0$ continuous creases by modifying the local polynomials~$\vec p(\vec \xi_j)$, and, as necessary, the chart domains~$\hat \Omega_j$. Creases are introduced along the edges of the control mesh prescribed by the user. The local polynomials~$\vec p(\vec \xi_j)$ are chosen to be $C^0$ continuous across the creases. The new chart domains require new maps, i.e.~$\vec \varPsi^\textrm{II}_j \colon \vec \eta \in \Box \mapsto \vec \xi_j \in \hat \Omega_j$, from the reference element to the chart domains. In turn, the blending functions~$w_j(\vec \xi_j)$ are combined from tensor-product B-splines using the new maps~$\vec \varPsi^\textrm{II}_j$. The remaining steps in manifold-based basis construction are identical to the non-creased case. 

\begin{figure}
  \centering
  \subfloat[][Type 1 crease]  
  { 
  	\includegraphics[scale=0.925]{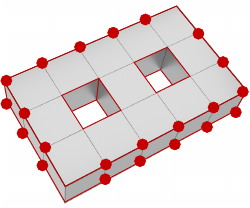}   \label{fig:blockTagA}
  }
  \hfil
  \subfloat[][Type 2 crease] 
  {
  	\includegraphics[scale=0.925]{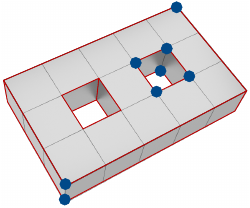}
  } 
  \hfil
  \subfloat[][Type 3 crease] 
  {
  	\includegraphics[scale=0.925]{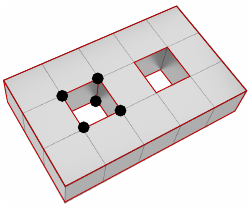}
  }
  \caption{Genus 2 surface with creased edges and sharp corners. The creased edges are marked in red. Each subfigure indicates the type of crease treatment to be used around the marked vertices.}
  \label{fig:blockTag}
\end{figure}

Depending on the number and arrangement of creases meeting at a vertex three different crease types are possible. As discussed in the following, each of them requires a different treatment. In Figure~\ref{fig:blockTag} the three different crease types are illustrated with the help of a sample geometry. 

%

\subsection{Crease Type 1: Rotationally symmetric chart domains}    
\label{sec:type1chart}                      
\begin{figure}
  \centering
  \subfloat[][$v_j=6$ and $k_j=2$] 
  {
  	\includegraphics[scale=0.875]{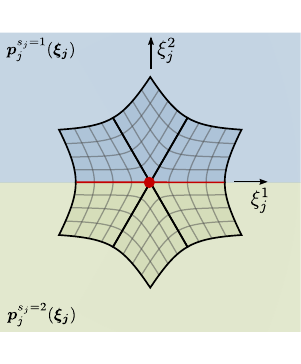}
  	\label{fig:v6c2chart}
  }
  \hfil
  \subfloat[][$v_j=6$ and $k_j=3$] 
  {
  	\includegraphics[scale=0.875]{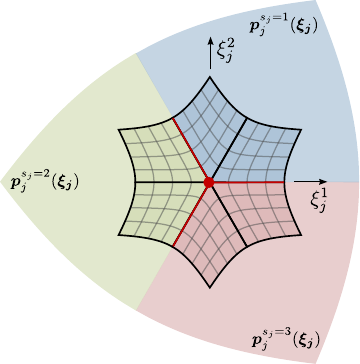}
  	\label{fig:v3c3chart}
  } 
  \caption{Two examples of Type 1 creases. The $k_j$ creased edges tagged by the user (marked in red) divide the chart domain into $k_j$ equiangular sectors shaded in different colours. Notice the rotational symmetry of the crease sectors. }
  \label{fig:chartType1}
\end{figure}

%
Type 1 creases have, as shown in Figure~\ref{fig:chartType1}, a rotationally symmetric chart domain with respect to the arrangement of the creased edges.  The centre vertex of a chart~$\hat \Omega_j$ can have $k_j \ge 2 $ attached edges tagged as crease\footnote{In this section, $k$ is used for the number of creased spoke edges, which is different from its meaning in $G^k$ and $C^k$ continuity.}. The~$k_j$ crease edges split the chart domain into~$k_j$ equiangular sectors, and the number of elements in each of the sectors must be the same. 

To introduce the basic idea in defining a suitable local approximant $f_j(\vec \xi_j)$, at first the case with two creased edges,~$k_j = 2$  , as depicted in Figure~\ref{fig:v6c2chart} is considered. The crease splits the chart domain into the two sectors~$s_j=1$ and~$s_j=2$ along the~$\xi_j^1$ coordinate axis. The local approximant is chosen as the piecewise continuous function
\begin{equation} \label{eq:symmetric}
\begin{split}
	& f_j(\vec \xi_j) = 
	\begin{cases}
		f^{s_j=1}_j(\vec \xi_j) & \text{ if } \xi^2_j \geq 0 \\
		f^{s_j=2}_j(\vec \xi_j) & \text{ if } \xi^2_j <0 \\
	\end{cases} 
	&\text{with} \quad f^{s_j=1}_j(\xi^1_j, \, \xi^2_j = 0) = f^{s_j=2}_j(\xi^1_j, \, \xi^2_j = 0)  \, .
\end{split}
\end{equation}
As in the non-creased case both approximants are given by  
\begin{equation}
	f^{s_j}_j (\vec \xi_j) = \vec p_j^{s_j}(\vec \xi_j) \cdot \vec \alpha_j^{s_j} \, . 
\end{equation}
The coefficients~$\vec \alpha_j^{s_j=1}$ and~$\vec \alpha_j^{s_j=2}$ must be matched so that~$ f_j(\vec \xi_j) $ is piecewise continuous.  Especially for a tensor-product Lagrange basis~$\vec p_j^{s_j} (\vec \xi_j)$, it is straightforward to determine a new set of coefficients~$\vec \alpha_j$ from the coefficients~$\vec \alpha_j^{s_j=1}$ and~$\vec \alpha_j^{s_j=2}$ which ensure piecewise continuity. The two Lagrange approximants must share the same coefficients along the common crease edge.

In chart domains with three or more creases,~$k_j \geq 3$, the bases~$\vec{p}^{s_j}_j(\vec{\xi}_j)$ are obtained by mapping a tensor-product Lagrange basis~$\vec{b}(\vec{\theta})$ to each sector. The respective local approximants $f^{s_j}_j(\vec{\xi}_j)$ in two adjacent sectors meet $C^0$~continuously because they share the same coefficients along the common crease edge. The Lagrange basis~$\vec b (\vec \theta)$ is mapped onto a sector~$s_j$ using the quasi-conformal map~\eqref{eq:conformalMap} introduced earlier 
\begin{equation} \label{eq:lagrangeToSector}
	\vec {p}^{s_j}_j(\vec \xi_j) = \vec b(\vec{\theta}) \quad \text{with}\quad \vec{\vec \theta} = \left(\vec{\Lambda_{qr}}^\textrm{I}\right)^{-1}(\vec{\xi}_j ; \, k_j, \, s_j)  \, .
\end{equation}
This mapping is best understood in conjunction with Figure~\ref{fig:v3c3chart}. Each of the three sectors in a different colour represents the image of a square mapped with~\mbox{$ \vec{\Lambda^\textrm{I}_{qr}} (\vec \theta; \, k_j, \, s_j) \colon \vec \theta \in \Box \mapsto \vec \xi_j \in \hat \Omega_j$} with $s \in \{ 1, \, 2, \, 3\}$. Obviously, the polynomial degree of~$\vec b(\vec \theta)$ must be chosen such that it is less than or equal to the number of vertices in each sector. More vertices have to be added, e.g. by refining the elements in each chart domain, if a higher degree basis is requested. Although we chose in our implementation a Lagrange basis for $\vec{b}(\vec \theta)$ it is possible to use other boundary interpolating tensor-product basis, like the Bernstein basis.

\subsection{Crease Type 2: Asymmetric chart domains \label{sec:type2chart}}                          
%
Type 2 creases have  a rotationally asymmetric  chart domain with respect to the arrangement of crease edges. Again, the centre vertex of a chart~$\hat \Omega_j$ can have~$k_j \ge 2 $ attached edges tagged as crease. However, the~$k_j$ crease edges split the chart domain into~$k_j$ sectors, which are not equiangular and have different number of elements. This makes it difficult to establish $C^0$ continuous local approximants~$f_j(\vec \xi_j) $ following the approach introduced for Type 1 creases. 

\begin{figure}
  \centering
  \includegraphics[scale=0.875]{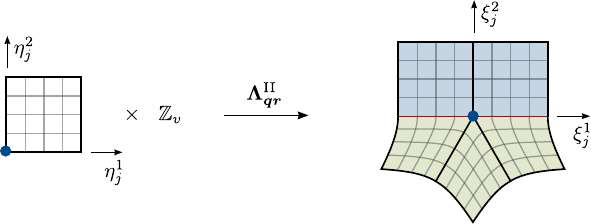}
  \caption{Example of  a Type 2 crease. The centre vertex has the valence $v_j=5$ and there are $k_j=2$ crease edges (marked in red). The number of elements in each of the two sectors is different, i.e. two and three respectively.}
  \label{fig:v5c2map}
\end{figure}

\begin{figure}
  \centering
  \includegraphics[scale=0.875]{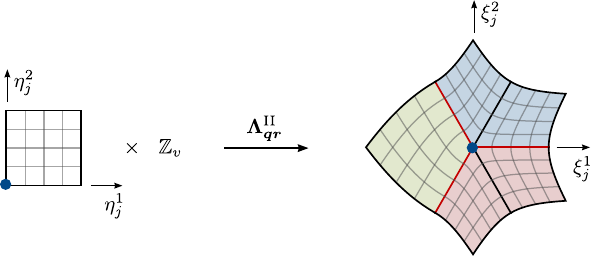}
  \caption{Example of a Type 2 crease. The centre vertex has the valence $v_j=5$ and there are $k_j=3$ crease edges (marked in red). The number of elements in each of the three sectors is different, i.e. two, one and two respectively.}
  \label{fig:v5c3map}
\end{figure}

\begin{figure}
  \centering
  \subfloat[][Sector number $s_j$] 
  {
  	\includegraphics[scale=0.875]{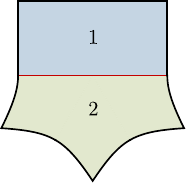}
  }
  \hfil
  \subfloat[][Number of elements in a sector $l_s$] 
  {
  	\includegraphics[scale=0.875]{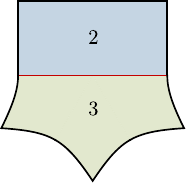}
  } 
  \hfil
  \subfloat[][Local element  number $m_s$] 
  {
  	\includegraphics[scale=0.875]{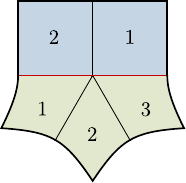}
  }   
  \caption{Illustration of the numbering used to define the map $\vec{\Lambda^\textrm{II}_{qr}}$ in~\eqref{eq:conformalMap2} with the help of the Type 2 crease chart shown in Figure~\ref{fig:v5c2map}. The number of creases $k_j=2$ takes the same value for every element.}
  \label{fig:v5c2numbering}
\end{figure}

Therefore, Type 2 creases are first mapped onto a chart domain, which is equiangular with respect to the arrangement of crease edges, see Figures~\ref{fig:v5c2map} and~\ref{fig:v5c3map}. This is possible in a manifold-based approach because the chart domains~$\hat \Omega_j$ can be freely chosen, subject to some smoothness and invertibility constraints. The respective maps \mbox{$\vec{\varPsi}^{\textrm{II}}_j  \colon \vec \eta \in \Box \mapsto \vec{\xi}_j \in \hat \Omega_j$}  have the same structure like for non-creased charts~\eqref{eq:composedMap}, i.e.,  
\begin{equation}\label{eq:composedMap2}
	\vec{\varPsi}^{\textrm{II}}_j = \vec{\Lambda_{qr}}^{\textrm{II}} \circ \vec{\Lambda_{tr}} \, , 
\end{equation}
where~$ \vec{\Lambda_{tr}}$  consists as defined in~\eqref{eq:linearMap} of  a translation and rotation and~$\vec{\Lambda_{qr}}^{\textrm{II}}$ is a quasi-conformal map with the arguments of~$ \vec{\Lambda^\textrm{I}_{qr}}$ in~\eqref{eq:conformalMap} replaced according to 
\begin{equation}
	\label{eq:conformalMap2}
	\vec{\Lambda_{qr}}^{\textrm{II}}(\vec{\eta}; \, k_j, \, s_j, \, l_s, \, m_s) = \vec{\Lambda^\textrm{I}_{qr}}(\vec{\eta}; \, k_jl_s, \, (s_j-1)l_s+m_s) \, .
\end{equation}
Here, $k_j$ is the number of creases, $s_j \in \{1, \, 2, \, \dotsc, \, k_j \}$ is the sector number, $l_s$ is the number of elements in the sector~$s_j$ and~$m_s \in \{1, \,  \dotsc, \,  l_s \}$ is the local element number in the sector~$s_j$.  The meaning of the introduced variables is further clarified in Figure~\ref{fig:v5c2numbering}. After the equiangular chart domain for Type 2 creases is established, as in Figures~\ref{fig:v5c2map} and~\ref{fig:v5c3map}, the local approximants~$f_j^{s_j}(\vec \xi_j)$ in each sector are obtained following the approach for a Type 1 crease chart. That is, they are defined as in~\eqref{eq:symmetric} in case of two creases or otherwise mapped according to~\eqref{eq:lagrangeToSector} from a tensor-product Lagrange basis~$\vec b(\vec \theta)$. The local approximant~$f_j(\vec \xi_j)$, composed from the local approximants in each of the sectors~$f_j^{s_j}(\vec \xi_j)$, is~$C^0$ continuous across the crease edges when the approximants share the same coefficients along the edges.

\subsection{Crease Type 3: Chart domains with concave corners \label{sec:type3chart}}                          
\begin{figure}
  \centering
  \includegraphics[scale=1.2]{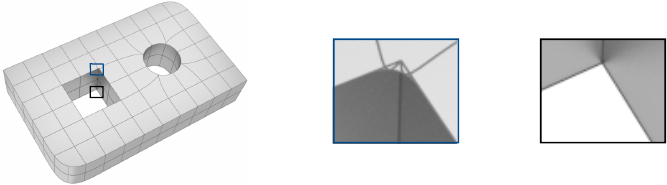}
  \caption{Behaviour of the manifold surface at concave corners. Edges are tagged as crease according to Figure~\ref{fig:blockTag}. The close-up image in the blue box (middle) shows  that the surface develops a fold when incorrectly a Type 2 mapping is used at the concave corner. The same behaviour can be observed when, for instance, Lagrange or B\'ezier elements are mapped to non-convex domains. The close-up image in the black box (right) shows that using a Type 3 mapping does not lead to a fold.}
  \label{fig:concaveDemo}
\end{figure}

%
The arrangement of creases in Type 2 chart domains can lead sometimes to non-convex sectors on the chart domain. However, there is no regular mapping from a convex to a non-convex domain with a non-zero determinant of the Jacobian. Therefore, when the map~$\vec{\varPsi}^{\textrm{II}}_j$ defined in~\eqref{eq:composedMap2} is used in a non-convex sector it leads to an unwanted folding-over of the surface as visible in Figure~\ref{fig:concaveDemo}.  This behaviour is similar to the one observed with isoparametrically mapped non-convex Lagrange finite elements. 

\begin{figure} 
  \centering
  \includegraphics[scale=0.875]{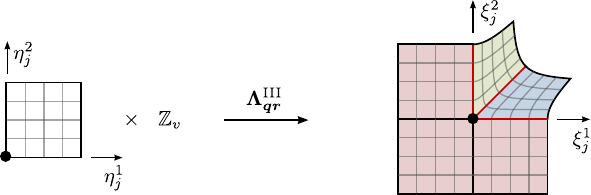}
  \caption{Example of Type 3 crease chart domain. The centre vertex has valence $v_j=5$ and there are $k_j=3$ creased edges (marked in red). One of the three sectors is non-convex (highlighted in red). }
  \label{fig:concavemap}
\end{figure}

\begin{figure}
  \centering
  \subfloat[][Sector number $s_j$] 
  {
  	\includegraphics[scale=0.875]{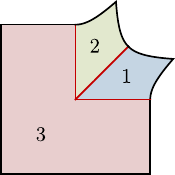}
  }
  \hfil
  \subfloat[][Number of elements in a sector $l_s$] 
  {
  	\includegraphics[scale=0.875]{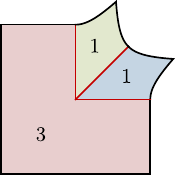}
  } 
  \hfil
  \subfloat[][Local element  number $m_s$] 
  {
  	\includegraphics[scale=0.875]{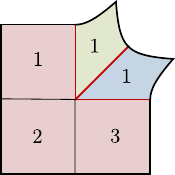}
  }   
  \caption{Illustration of the numbering used to define the map $\vec{\varPsi_c}^{\textrm{III}}$ in~\eqref{eq:concavemap} taking the example of the Type 3 crease chart shown in Figure~\ref{fig:concavemap}. The number of creases $k_j=3$ takes the same value for every element.}
  \label{fig:v5c3numbering}
\end{figure}

The  folding-over can be avoided by composing the non-convex sector from several smoothly connected convex pieces. To this end, the chart domain depicted in Figure~\ref{fig:concavemap} is introduced. The respective map~\mbox{$\vec{\varPsi}^{\textrm{III}}_j  \colon \vec \eta \in \Box \mapsto \vec{\xi}_j \in \hat \Omega_j$} has the same structure like the previously introduced maps, c.f.~\eqref{eq:composedMap} and~\eqref{eq:composedMap2}, 
\begin{equation}\label{eq:composedMap3}
	\vec{\varPsi}^{\textrm{III}}_j = \vec{\Lambda_{qr}}^{\textrm{III}} \circ \vec{\Lambda_{tr}}  
\end{equation}
with the quasi-conformal map~$\vec{\Lambda_{qr}}^{\textrm{III}}$ obtained by replacing the arguments of~$ \vec{\Lambda^\textrm{I}_{qr}}$ in~\eqref{eq:conformalMap} in the following way
\begin{equation} \label{eq:concavemap}
\begin{split}
	& \vec{\Lambda_{qr}}^{\textrm{III}} (\vec{\eta}; \, k_j, \, s_j, \, l_s, \, m_s) = 
	\begin{cases}
		\vec{\Lambda^\textrm{I}_{qr}}(\vec{\eta}; \, 4(k_j-1)l_s, \,  (s_j-1)l_s+m_s) & \text{ if } s_j < k_j  \\
		\vec{\Lambda^\textrm{I}_{qr}}(\vec{\eta}; \, 4l_s/3, m_s) \cdot e^{i \pi /2} & \text{ if } s_j = k_j \\
	\end{cases} \, , 
\end{split}
\end{equation}
where $k_j$ is the number of creases, $s_j \in \{1, \, 2, \, \dotsc , \, k_j \}$ is the sector number, $l_s$ is the number of elements in a sector and  $m_s$ is the local element number in a sector. All the mentioned variables have the same meaning as in Section~\ref{sec:type2chart} and are clarified in Figure~\ref{fig:v5c3numbering}. As shown in Figure~\ref{fig:concavemap}, with the introduced map~$\vec{\varPsi}^{\textrm{III}}_j$ the concave sector with $s_j = k_j$ is mapped to an L-shaped domain and all other sectors are mapped to the first quadrant. For obtaining the local approximant~$f_j(\vec \xi_j)$ a similar approach as discussed for the Type 1 and 2 crease charts is followed  with the exception of the non-convex sector. The local approximants~$f_j^{s_j<k_j}(\vec \xi_j)$ on each convex sector are mapped according to~\eqref{eq:lagrangeToSector}  from a tensor-product Lagrange basis~$\vec b (\vec \theta)$ using~$\vec \theta =  (\Lambda^{I}_{qr} )^{-1} (\vec \xi_j;  \, 4(k_j-1), \, s_j)$. The local approximant on the non-convex sector~$f_j^{s_j=k_j}(\vec \xi_j)$ is composed out of three (cubic) B\'ezier patches that are smoothly connected within the sector while they are matched $C^0$ continuously with polynomial pieces in other sectors across the crease edges.

\subsection{Examples of manifold surfaces with creases}
\begin{figure}
	\centering 
	\includegraphics[scale=1.075]{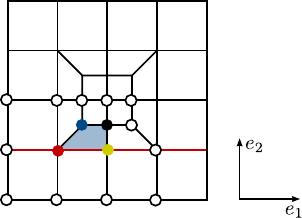}	
	\caption{An unstructured mesh with extraordinary vertices of valence $v \in \{3, \,  5 \} $. The shaded element is overlapped by four charts and the union of the four charts have 16 unique vertices in total. Edges tagged as crease are marked in red. \label{fig:mesh2}}
	\vspace*{5\floatsep}
%
	\centering 
	\includegraphics[scale=0.875]{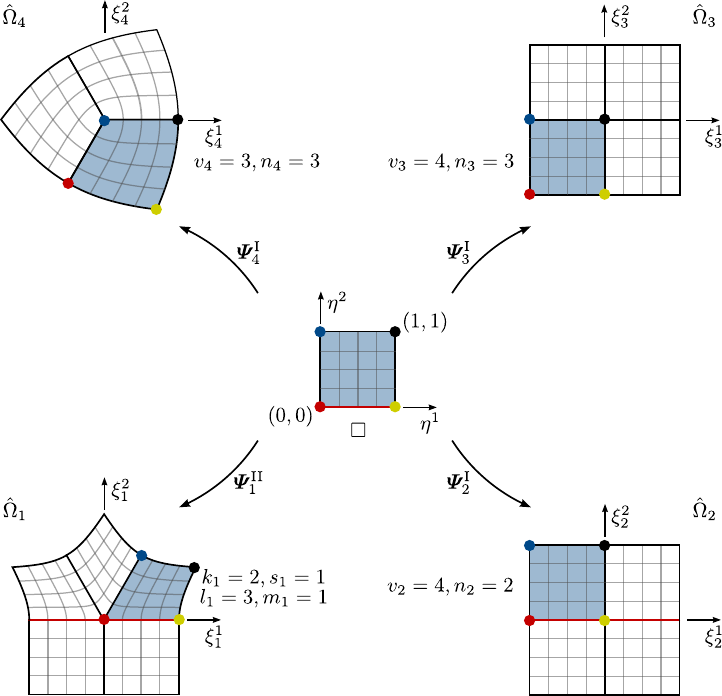}	
	\caption{The reference element $\Box$ (centre) and its four overlapping chart domain $\hat \Omega_j$ with $j\in \{1, \, 2, \, 3, \,4\}$. The reference element is mapped with $\vec \varPsi_j$ to the blue shaded elements in the four charts $\hat \Omega_j$. Notice that the parameter lines in the charts (with $\eta^1 = \text {const.}$ and $\eta^2= \text{const.}$ ) are always orthogonal to the spoke edges which guarantees that the derivatives of the maps $\vec \varPsi_j$ are  continuous over the entire chart. \label{fig:manifoldPUM2}}
\end{figure}

%
We consider the construction of the manifold-based basis functions for one of the elements in the unstructured quadrilateral control mesh shown in Figure~\ref{fig:mesh2}. The element is adjacent to a prescribed crease and belongs to four overlapping chart domains as illustrated in Figure~\ref{fig:manifoldPUM2}. Two of the chart domains,~$\hat \Omega_1$ and~$\hat \Omega_2$,  contain each two crease edges and the other two,~$\hat \Omega_3$ and~$\hat \Omega_4$, contain no crease edges. In the chart domain~$\hat \Omega_2$ the arrangement of the crease edges is rotationally symmetric so that it is Type~1 and in the chart domain~$\hat \Omega_1$ it is asymmetric so that it is Type~2.  For the chart~$\hat \Omega_1$ the mapping~\eqref{eq:composedMap2} and for the chart~$\hat \Omega_2$ the mapping~\eqref{eq:composedMap} is to be used. Hence, in comparison to the smooth case illustrated in Figure~\ref{fig:manifoldPUM} only the mapping for the chart~$\hat \Omega_1$ is different while the other mappings are the same. The presence of creased edges only changes the chart parametrisation for the construction of local basis $\vec{p}_j(\vec{\xi}_j)$. The construction of blending functions $w_j(\vec \xi_j)$ is the same for all charts regardless of the presence of crease edges. The tensor-product B-spline blending functions are defined in the reference element and mapped to each element in the chart domains.  

\begin{figure}
  \centering
  \includegraphics[scale=0.9]{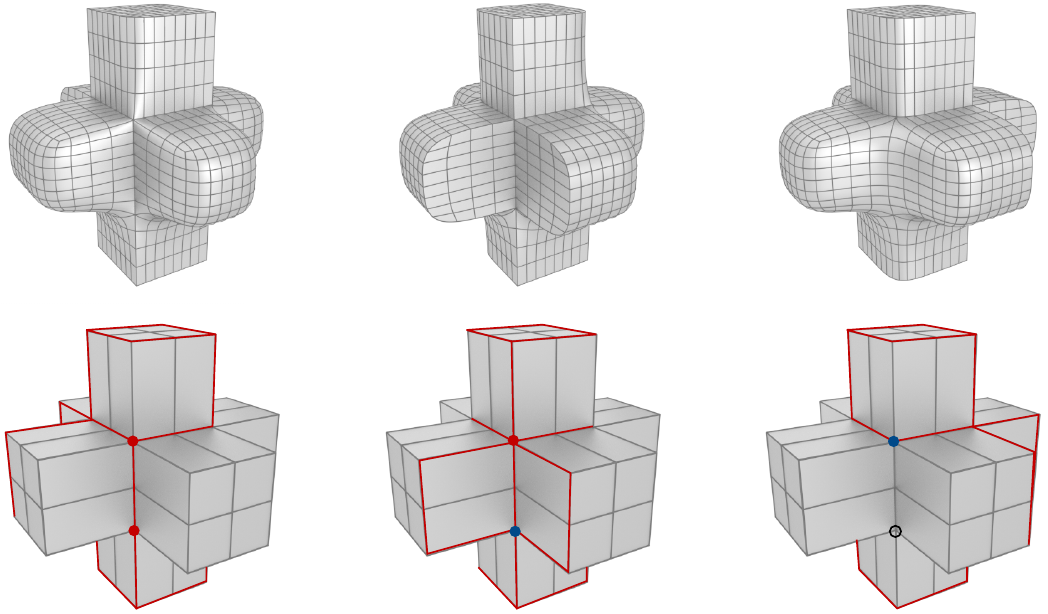}
  \caption{Three different views of a manifold surface (top) approximating a given control mesh (bottom). In the control mesh, crease edges are marked in red. The mesh includes eight valence~$v=6$ vertices and the corresponding charts are categorised as a smooth chart  (empty dot) and Type~1 (red dot) or Type~2 (blue dot) crease charts. Across the indicated crease edges the manifold surface is $C^0$ continuous and otherwise it is $C^2$ continuous.}
  \label{fig:geoExample}
\end{figure}

In Figure~\ref{fig:geoExample} the application of the manifold-based basis functions in representing a smooth surface with sharp features is demonstrated. The blending functions are assembled from tensor-product cubic B-splines and the local polynomial  basis is a  cubic tensor-product Lagrange polynomial. The relatively coarse unstructured quadrilateral mesh has a number of tagged crease edges. The location and arrangement of creases have been chosen in order to obtain as many as possible distinct chart configurations.  The eight charts with valence $v=6$ include a smooth chart, four Type 1 crease charts and three Type 2 crease charts. The rendered manifold surface~$\vec x(\vec \eta)$ is obtained by  first projecting the control vertices~$\vec x_I$ on a Catmull-Clark subdivision surface 
\begin{equation}
	\vec x^S_I =\sum_I   L_{IJ} \vec x_J   \, ,
\end{equation}	 
where $L_{IJ}$ is the limit matrix (or, mask) of the subdivision surface. The multiplication of the control vertices with the limit matrix involves only the one-neighbourhoods of the vertices~$\vec x_I$ and projects them onto a Catmull-Clark limit surface. The entries of the limit matrix depend on the valence of the vertex and can be found, for instance, in~\cite{Biermann:2000aa}.  There are also so-called tuned subdivision masks available which can provide better shapes~\cite{zhang2018subdivision,augsdorfer2006tuning,Peters:2008aa}. Moreover, note that the Catmull-Clark subdivision surface reduces to cubic B-splines on structured meshes. With the projected control vertices the image of a reference element~$\Box$ on the manifold surface is obtained as
\begin{equation}
	\vec x(\vec \eta) = \sum_I   N_I (\vec \eta)   \vec x_I^S    \, ,
\end{equation}	
where~$N_I (\vec \eta)$ are the introduced manifold-based basis functions.  The manifold construction ensures that the images of the reference elements on the surface are connected with the required smoothness. The projection of the control vertices onto a Catmull-Clark surface is only used to obtain a visually pleasing surface.  It is inconsequential when manifold-based basis functions are used as finite element basis functions. 

Finally, evidently, the geometry shown in Figure 18 can be modelled using subdivision surfaces, or even using a collection of B-spline patches, without manifold constructions. The obtained smooth manifold-based basis functions are however essential for the finite element discretisation, in particular, of higher-order partial differential equations, such as the Kirchhoff-Love thin shell equations. 

%% file: mechanics.tex
\section{Thin-shell formulation and discretisation} \label{sec:thinshell}         
%
This section provides a brief review of the Kirchhoff-Love thin-shell equations and their discretisation with manifold-based basis functions. Also, the boundary and crease constraints used in the presented numerical examples are introduced. Further details of our specific Kirchhoff-Love implementation may be found in~\cite{Cirak:2001aa, Cirak:2011aa, Ciarlet:2005aa}. 

We denote the reference and deformed mid-surfaces of the shell with~$\Omega_0$ and~$\Omega$, respectively, and the corresponding boundaries with~$\Gamma_0=\partial \Omega_0$ and~$\Gamma = \partial \Omega$. Any configuration of the shell is assumed to be defined as 
\begin{equation}
	\vec{\varphi}(\eta^1, \, \eta^2, \, \eta^3) = \vec{x}(\eta^1, \, \eta^2) + \eta^3 \vec{a}_3(\eta^1, \, \eta^2)  \quad \text{ with} \quad  \eta^3 \in \left [-\frac{t}{2}, \, \frac{t}{2}  \right ]\,,
\end{equation}
where~$\vec{\varphi}(\eta^1, \,  \eta^2, \, \eta^3)$  is the position vector of a material point with the convective coordinates~$(\eta^1, \, \eta^2, \, \eta^3)$. Similarly,~$\vec{x}(\eta^1, \, \eta^2)$ is the position vector of a material point with the convective coordinates~$(\eta^1, \, \eta^2)$ on the shell mid-surface. Furthermore,~$\vec a_3(\eta^1, \, \eta^2) $ is the unit normal to the mid-surface~$\vec{x}(\eta^1, \, \eta^2)$  and~$t$ is the shell thickness. The corresponding vectors in the reference configuration are denoted with uppercase letters~$\vec{\varPhi}(\eta^1, \, \eta^2,  \, \eta^3)$, $\vec X(\eta^1, \,  \eta^2)$ and~$ \vec A_3(\eta^1, \, \eta^2)$ so that the reference configuration reads
\begin{equation}
	\vec{\varPhi}(\eta^1, \, \eta^2, \, \eta^3) = \vec{X}(\eta^1, \, \eta^2) + \eta^3 \vec{A}_3(\eta^1, \, \eta^2) \,.
\end{equation}
The covariant basis vectors of the tangent space of the mid-surface are defined as\footnote{Throughout Section~\ref{sec:thinshell}, the Greek indices take the values~$\{1, \, 2 \}$, the lowercase Latin indices take the values~$\{1, \, 2, \, 3 \}$, a comma denotes differentiation and the summation convention over repeated indices applies.}
\begin{equation}
	\vec{A}_\alpha = \frac{\partial \vec{X}}{\partial \eta^\alpha} \quad \text{and}  \quad \vec{a}_\alpha = \frac{\partial \vec{x}}{\partial \eta^\alpha} 
\end{equation}
with the corresponding normal vectors 
\begin{equation}
	\vec{A}_3 = \frac{\vec{A}_1 \times \vec{A}_2}{|\vec{A}_1 \times \vec{A}_2|} \quad \text{and} \quad \vec{a}_3 = \frac{\vec{a}_1 \times \vec{a}_2}{|\vec{a}_1 \times \vec{a}_2|} \, .
\end{equation}
The contravariant basis vectors~$\vec A^\alpha$ and~$\vec a^\alpha$  follow from the relations
\begin{equation}
	\vec A^\alpha \cdot \vec A_\beta = \delta_\alpha^\beta \quad \text{and} \quad \vec a^\alpha \cdot \vec a_\beta = \delta_\alpha^\beta \, ,
\end{equation}
where~$\delta_\alpha^\beta$ is the Kronecker delta. The contravariant metric of the reference configuration needed in the following is defined as 
\begin{equation}
	\vec A^{\alpha \beta} = \vec A^\alpha \cdot \vec A^\beta \, .
\end{equation}

The deformation gradient can now be expressed as
\begin{equation}
	\vec{F} = \frac{\partial \vec{\varphi}}{\partial {\vec{\varPhi}}} = \frac{\partial \vec{\varphi}}{\partial \eta^j}  \otimes \frac{\partial \eta^j}{\partial \vec{\varPhi}} =   \vec a_j \otimes \vec A^j \,,
\end{equation}
and a straightforward calculation yields the Green-Lagrange strain tensor 
\begin{equation}
	\vec{E} = \frac{1}{2}(\vec{F}^\trans \vec{F} - \vec{I})  = \vec{\alpha} + \eta^3\vec{\beta} + \left (\eta^3 \right )^2 \dotsc 
\end{equation}
with the membrane and bending strain tensors 
\begin{subequations}
\begin{align}
	\vec{\alpha} &= \frac{1}{2}(\vec{a}_\alpha \cdot \vec{a}_\beta - \vec{A}_\alpha \cdot \vec{A}_\beta ) \, \vec{A}^\alpha \otimes \vec{A}^\beta  \\
	\vec{\beta} &= \frac{1}{2}(\vec{a}_\alpha \cdot \vec{a}_{3,\beta} + \vec{a}_\beta \cdot \vec{a}_{3,\alpha}- \vec{A}_\alpha \cdot \vec{A}_{3,\beta} - \vec{A}_\beta \cdot \vec{A}_{3,\alpha}) \, \vec{A}^\alpha \otimes \vec{A}^\beta \,.
\end{align}
\end{subequations}
As usual, the quadratic terms in the Green-Lagrange strain tensor~$\vec E$ are neglected for thin shells. 

The potential energy of a hyper-elastic shell is given by 
\begin{align} \label{eq:energy}
	\begin{split}
	\Pi(\vec x) & = \int_{\Omega_0} W(\vec \alpha, \, \vec \beta) \D \Omega_0 + \Pi^{\text{ext}} (\vec{x})   \\
	& =  \int_{\Omega_0}  \left ( \frac{1}{2} \frac{Et}{1-\nu^2}H^{\alpha \beta \gamma \delta}\alpha_{\alpha\beta} \alpha_{\gamma\delta} + \frac{1}{2} \frac{Et^3}{12(1-\nu^2)}H^{\alpha \beta \gamma \delta}\beta_{\alpha\beta} \beta_{\gamma\delta} \right )  \D \Omega_0 + \Pi^{\text{ext}} (\vec{x}) \,,
	\end{split}
\end{align}
where~$W(\vec \alpha, \, \vec \beta) $ is the internal energy density and~$\Pi^{\text{ext}}(\vec{x}) $ is the potential of the external forces. The internal energy density depends, in addition to the two strain tensors~$\vec \alpha$ and~$\vec \beta$, on the Young's modulus~$E$, the Poisson's ratio~$\nu$ and a mainly geometry related fourth-order tensor with the components
\begin{equation*}
	H^{\alpha \beta \gamma \delta} = \nu A^{\alpha\beta} A^{\gamma\delta} + \frac{1}{2}(1-v)(A^{\alpha\gamma}A^{\beta\delta}  + A^{\alpha\delta}A^{\beta\gamma}) \,.
\end{equation*}

When rigid joints or clamped boundaries are present, it is necessary to constrain the surface normals. The introduced manifold-based basis functions are interpolating at the boundaries so that it is straightforward to impose displacement boundary conditions. Along the clamped boundaries~$\Gamma_{0,c}$ the change of the mid-surface normal~$(\vec a_3 - \vec A_3)$ is required to be zero. Along the joints~$\Gamma_{0,j}$, the~$C^0$ continuity of the basis functions ensures that the displacements to the left~$(l)$ and to the right~$(r)$ are the same. When the joints are structurally rigid the relative change of the mid-surface normals~$\vec{a}_3^{(l)} \cdot \vec{a}_3^{(r)} - \vec{A}_3^{(l)} \cdot \vec{A}_3^{(r)}$  is constrained to be zero.  In our current implementation the rigid joint and clamped boundary constraints are  enforced with the penalty method 
\begin{equation}\label{eq:penalty}
	\Pi^{C}(\vec{x}) = \Pi(\vec{x}) + \underbrace{ \frac{\gamma_1}{2} \int_{\Gamma_{\text{r}}}  \left (\vec{a}_3^{(l)} \cdot \vec{a}_3^{(r)} - \vec{A}_3^{(l)} \cdot \vec{A}_3^{(r)} \right )^2 \D \Gamma}_\text{rigid joints}  + \underbrace{\frac{\gamma_2}{2} \int_{\Gamma_{\text{c}}}  \left (\vec{a}_3 - \vec{A}_3 \right )\cdot \left (\vec{a}_3 - \vec{A}_3 \right) \D \Gamma}_\text{clamped edges} \,,
\end{equation}
where $\gamma_1$ and $\gamma_2$ are penalty parameters. It is straightforward to consider alternative constrain enforcement techniques, like the Nitsche method~\cite{Embar:2010aa, guo2017parameter} or Lagrange multipliers~\cite{Cirak:2011aa, duong2017new}.

As usual, the discrete finite element equilibrium equations are derived by first writing the potential~\eqref{eq:penalty} as a sum over the set of  reference elements~$\{ \Box_e \}$ using the Jacobian~$ | \partial \vec X/\partial \vec \eta  |$. Here, the index~$e$ denotes the number of an element.  It is emphasised that each element corresponds to a quadrilateral in the control mesh, as in isogeometric analysis with B-splines. The charts considered in the manifold construction are only used for obtaining the basis functions. With the determined basis functions the reference and deformed mid-surface of a reference element~$\Box_e$ are approximated by 
\begin{equation}
	\vec{X}^h(\eta^1, \,  \eta^2) = \sum_{I}N_I (\eta^1, \,  \eta^2)\vec{X}_I   \quad  \text{and} \quad \vec{x}^h(\eta^1, \, \eta^2) = \sum_{I}N_I (\eta^1,  \, \eta^2)\vec{x}_I  \, ,
\end{equation}
where ~$\vec X_I$ and~$\vec  x_I$ are the coordinates of the control vertex coordinates. After introducing~$\vec{X}^h(\eta^1,  \, \eta^2) $ and~$\vec{x}^h(\eta^1,  \, \eta^2)$ into the element-wise expressed potential energy~\eqref{eq:penalty} and numerical integration the discrete equilibrium equations follow from the stationarity principle
\begin{equation}
	\frac{\partial \Pi(\vec{x}^h)}{\partial \vec{x}_I} = 	\frac{\partial \Pi^\text{int}(\vec{x}^h)}{\partial \vec{x}_I} + 	\frac{\partial \Pi^\text{ext}(\vec{x}^h)}{\partial \vec{x}_I} =  \vec{0} \, .
\end{equation}
We solve this set of nonlinear equations iteratively with the Newton-Raphson method. That is, at each iteration step~$(n)$ the linear equation 
\begin{equation}\label{eq:equationsystem}
	\frac{\partial^2 \Pi \left (\vec{x}_I^{(n-1)} \right )}{\partial \vec{x}_I \partial \vec{x}_J} \left (\vec{x}_J^{(n)} - \vec{x}_J^{(n-1)} \right )= - \frac{\partial \Pi \left (\vec{x}_I^{(n-1)} \right )}{\partial \vec{x}_I} 
\end{equation}
is solved to determine~$\vec x_J^{(n)}$ for a given~$\vec x_J^{(n-1)}$.   The expression $\partial^2 \Pi/\partial \vec{x}_I \partial \vec{x}_J$ on the left is the tangent stiffness matrix whereas the expression $\partial \Pi/\partial \vec{x}_I$ on the right is the residual vector. The full expression for~\eqref{eq:equationsystem} can be found in~\cite{Cirak:2000aa, Long:2012aa}.
%

%% file: examples.tex
\section{Examples}                          
%
We consider beam, plate and shell examples to demonstrate the utility and convergence of the proposed creased manifold-based basis functions in finite element analysis. In all the examples the  Kirchhoff-Love model reviewed in Section~\ref{sec:thinshell} is discretised with manifold-based basis functions. The blending functions are assembled from tensor-product cubic B-splines and the local polynomial basis functions are quadratic tensor-product polynomials.  In convergence studies all the finite element integrals are integrated with $9\times 9$ Gauss quadrature points. The used standard Gauss integration rules do not take into account that the introduced polynomial manifold-based basis functions consist in each element out of~$4 \times 4$ piecewise polynomials. For the geometrically nonlinear simulation of a pinched tube we use $3\times 3$ quadrature points to save computing time.

In constructing the basis functions all element edges on the domain boundaries are tagged as crease. This eliminates the need for ghost cells as originally used in~\cite{majeedCirak:2016}. As discussed, the number of basis functions in a chart domain has to be equal or less than the number of vertices. Hence, to be able to use a quadratic tensor-product basis in all charts, mid-face and mid-edge vertices are introduced. This increases the number of vertices in a chart domain from~$2 v_j+1$ to~$6 v_j+1$, where~$v_j$ is the valence of the centre vertex.  In all examples, the rigid joint and clamped boundary constraints  are enforced with the penalty method. The penalty parameters are chosen to be $\gamma = \mathcal{O}(E)$, where~$E$ is the Young's modulus. 

\subsection{Propped cantilever beam\label{sec:beam}}
\begin{figure}[]
	\centering 
	\includegraphics[scale=0.85]{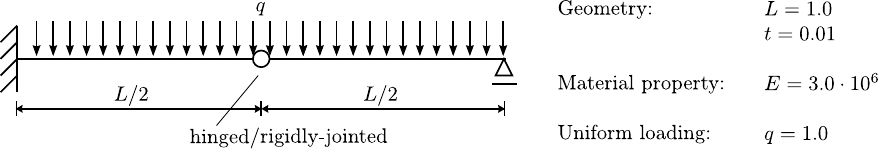}
	\vspace{0.01\textheight}	
	\caption{Geometry and loading of the propped cantilever beam. The left end is clamped and the right end is simply supported. Different treatments are considered for the hinge at the midspan. }
	\label{fig:beam}
\end{figure}

%
As an introductory example, we compute the deflection of a beam subjected to a uniformly distributed transversal loading, see Figure~\ref{fig:beam}. To be able to use bivariate basis functions the beam is modelled as a plate with uniform width and is discretised with a structured quadrilateral mesh. The chosen boundary conditions, loading and the Poisson ratio of~$\nu=0$ ensure that the structure deflects like a beam. Both the left and right ends of the beam are tagged as crease, eliminating the need for ghost cells. The left end is clamped so that the displacement and the change of the beam axis normal are required to be zero. The right end is simply supported so that only the displacement is required to be zero. 

In the computations the beam midspan is modelled either as continuous, hinged or rigidly-jointed hinged. The considered meshes have always a vertex at the midspan. For the continuous beam (with no hinge) the chart domain corresponding to the midspan vertex is treated in the same way like the other chart domains.  In the two cases with a hinge the midspan vertex  is tagged as a crease. In geometric terms, across a hinge the displacement is continuous but its derivative may be discontinuous. And, in structural terms, across a hinge only shear forces may be transferred while the bending moment at the hinge is  zero. In the case of the rigidly-jointed hinged beam, the two beam axis normals across the hinge are enforced to be weakly continuous (with the penalty method) so that  it becomes possible to transfer moments. 

\begin{figure}
	\centering 
  	\subfloat[][$L_2$ norm] 
  	{
  		\includegraphics[scale=0.45]{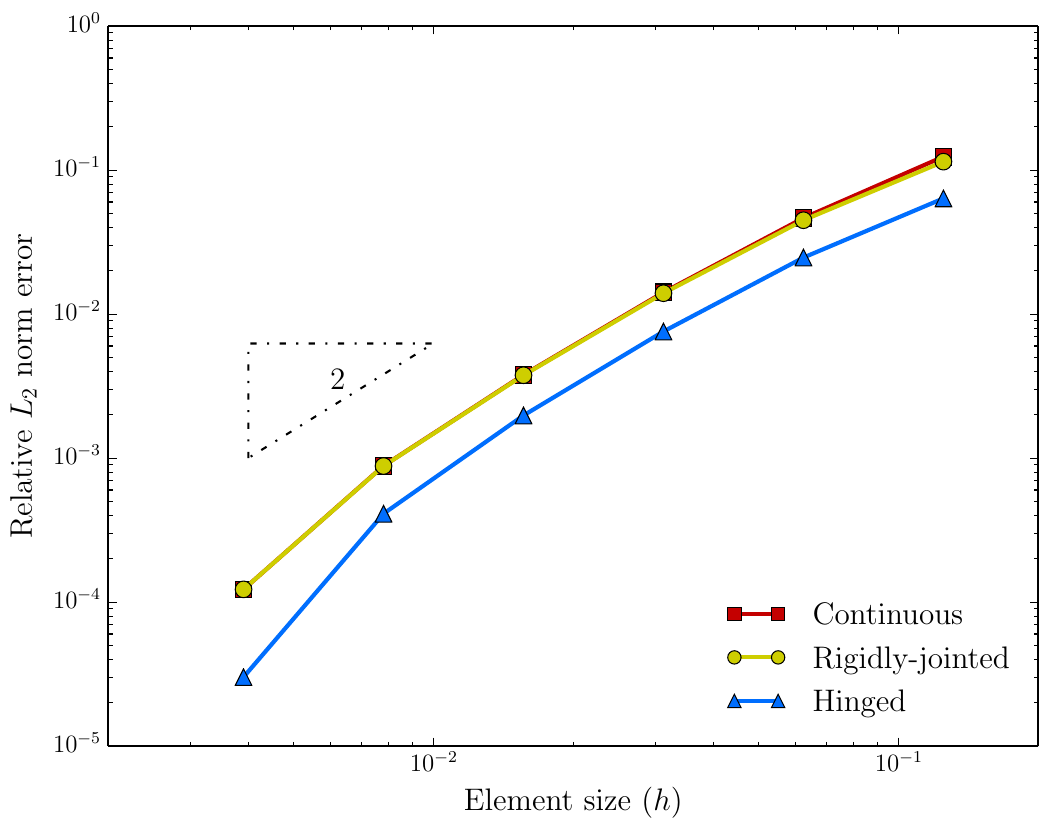}
  		\label{fig:beamErrorL2}
  	}
  	\subfloat[][Energy norm] 
  	{
  		\includegraphics[scale=0.45]{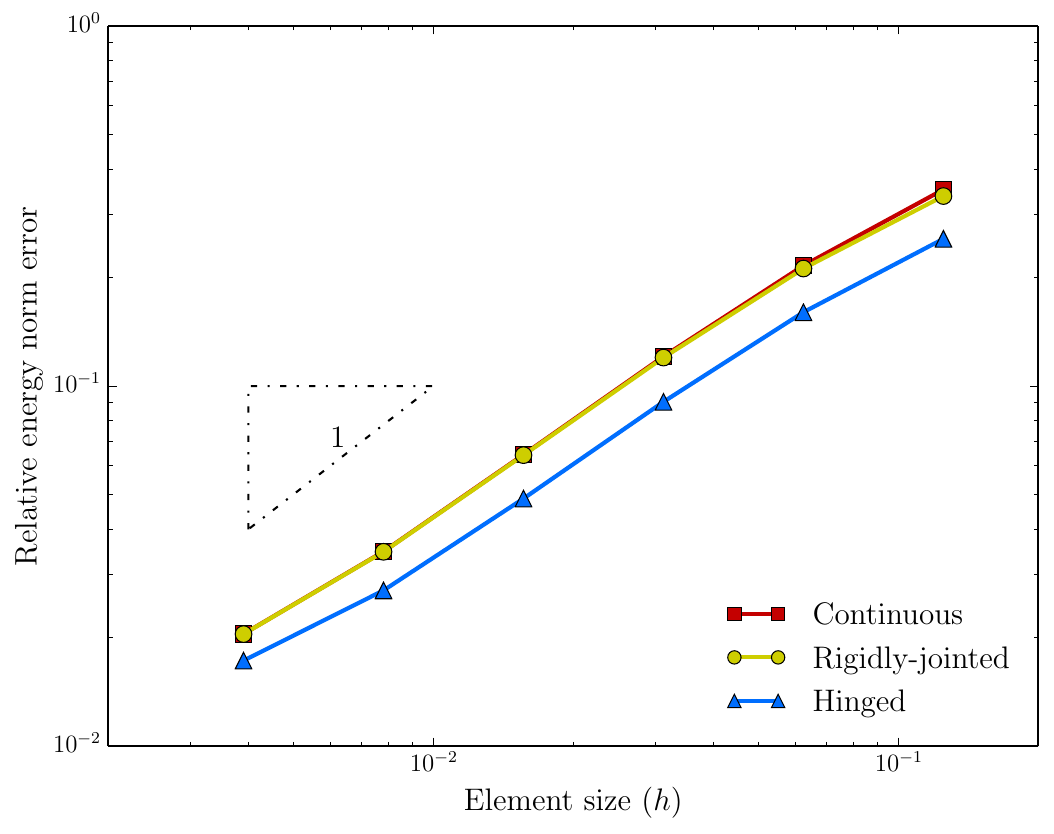}
  		\label{fig:beamErrorE}
 	} 
	\caption{Convergence of the relative $L_2$ and energy norm errors for the propped cantilever beam. The rigidly-jointed and hinged beams contain a hinge at the midspan while the continuous beam does not. For the rigidly-jointed beam the beam axis normals across the hinge are constrained to be weakly continuous. }
	\label{fig:beamError}
\end{figure}

As can be obtained by a straightforward integration the analytical deflection of the continuous propped cantilever is
\begin{equation} \label{eq:continuousBeam}
	u(x) = \frac{q}{4 Et^3}(2x^4 - 5Lx^3 + 3L^2x^2) \,,
\end{equation}
and the deflection of the propped cantilever with a hinge is  
\begin{equation} \label{eq:hingedBeam}
\begin{split}
	& u(x) = 
	\begin{cases}
		\dfrac{q}{2 Et^3}(x^4 - 3Lx^3 + 3L^2x^2) & \text{ if } x < L/2 \\
		\dfrac{q}{2 Et^3}(x^4 - 3Lx^3 + 3L^2x^2 -2L^3x + L^4) & \text{ if } x \geq L/2 \\
	\end{cases} \, . 
\end{split}
\end{equation}
In Figure~\ref{fig:beamError} the convergence of the relative $L_2$ norm and energy norm errors with decreasing mesh size are plotted. As can be seen, both $L_2$ and energy norm errors converge in all cases with the optimal rate. The errors for the beam with the hinge are slightly smaller. As to be expected, the continuous beam and rigidly-jointed hinged beam give very similar results. This confirms the soundness of the enforcement of rigid joint constraints with the penalty method. 

\subsection{Square plate  \label{sec:plate}}
\begin{figure}[]
	\centering 
	\includegraphics[scale=1]{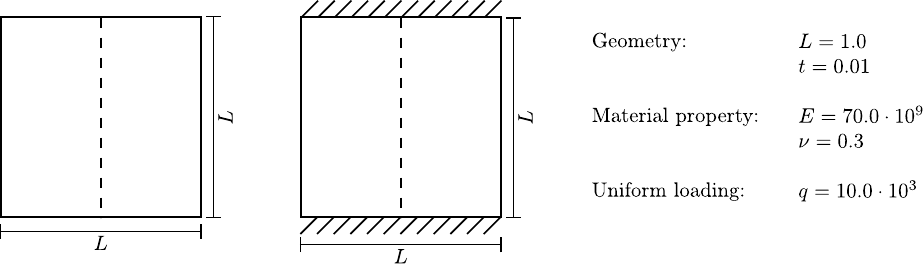}
	\caption{Definition of the square plate problem. Two types of boundary conditions are considered, with either all edges simply supported (left), or two opposite edges are simply supported while the other two are clamped (centre). The dashed vertical line at the centre indicates a hinge, which is present in some of the computations. If a hinge is present, the mid-surface normals across the hinge line are constrained to be weakly continuous.}
	\label{fig:plateLoading}
\end{figure}

%
Next, we   consider the deflection of a square plate with different boundary conditions subjected to uniformly distributed transversal loading, see Figure~\ref{fig:plateLoading}. Two different types of boundary conditions are considered. In the first case all boundary edges are simply supported and in the second case two opposite boundary edges are simply supported and the other two are clamped. At the simply supported boundary edges the displacements are constrained to be zero and at the clamped edges both the displacements and their derivatives are constrained to be zero. Similar to the propped cantilever example also a hinge line is introduced along the centre of the plate, which is considered to be present in only some of the computations. If a hinge line is considered to be present, it is always modelled as rigidly-jointed. 

\begin{figure}[]
	\centering 
	\subfloat[][Structured] {
	\includegraphics[scale=1.1]{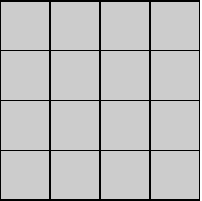}
	}
	\hfil
	\subfloat[][Unstructured] {
	\includegraphics[scale=1.1]{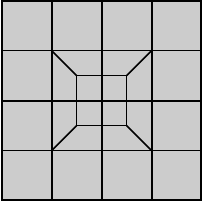}
	}
	\caption{Coarse control meshes used in square plate computations.}
	\label{fig:plateMesh}
\end{figure}

\begin{figure}[]
	\centering 
	\subfloat[][Simply supported] {
	\includegraphics[scale=0.25]{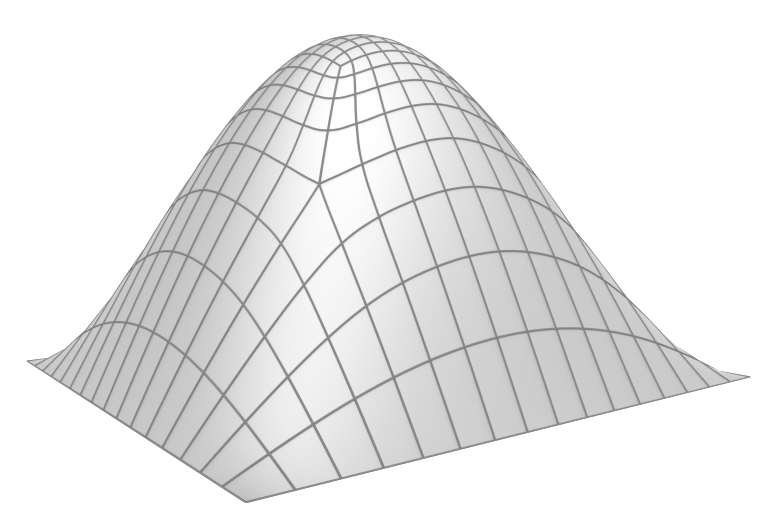}
	}
	\hfil
	\subfloat[][Partially clamped] {
	\includegraphics[scale=0.25]{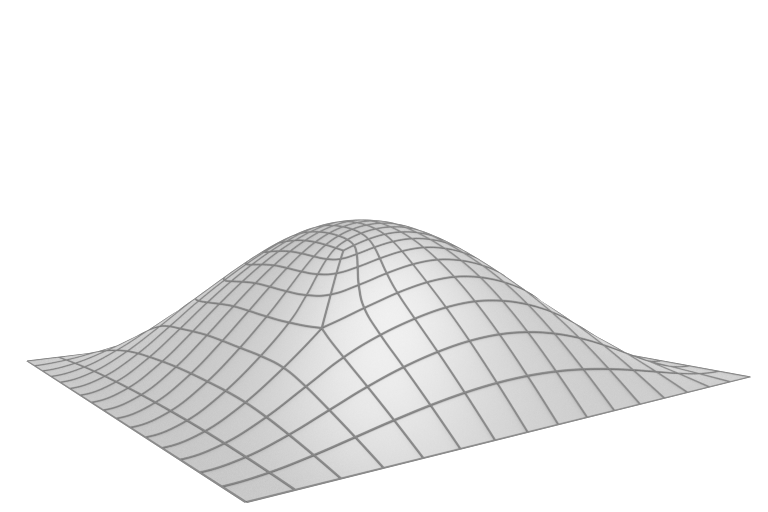}
	}	
	\caption{Deflected shapes of the square plate with two different boundary conditions. The deflections are scaled with the same factor in both cases.}
	\label{fig:plateDeformed}
\end{figure}

The structured and unstructured coarse meshes used in the computations are depicted in Figure~\ref{fig:plateMesh}. The unstructured mesh has eight extraordinary vertices, with four vertices of valence~$v_j=3$ and the other four of valence~$v_j=5$. In the convergence studies the refined meshes are obtained by subdividing the shown coarse meshes with the Catmull-Clark subdivision.  The analytical deflections for both boundary conditions can be found in Timoshenko and Woinowsky-Krieger~\cite[Chapters 5 and 6]{Timoshenko1959}. Two representative finite element solutions for the considered two boundary conditions are shown in Figure~\ref{fig:plateDeformed}. As to be expected the deflection of the plate with two clamped edges is much smaller than the deflection of the fully simply supported plate. 

\begin{figure}[]
	\centering 
  	\subfloat[][Structured mesh] 
  	{
  		\includegraphics[scale=0.45]{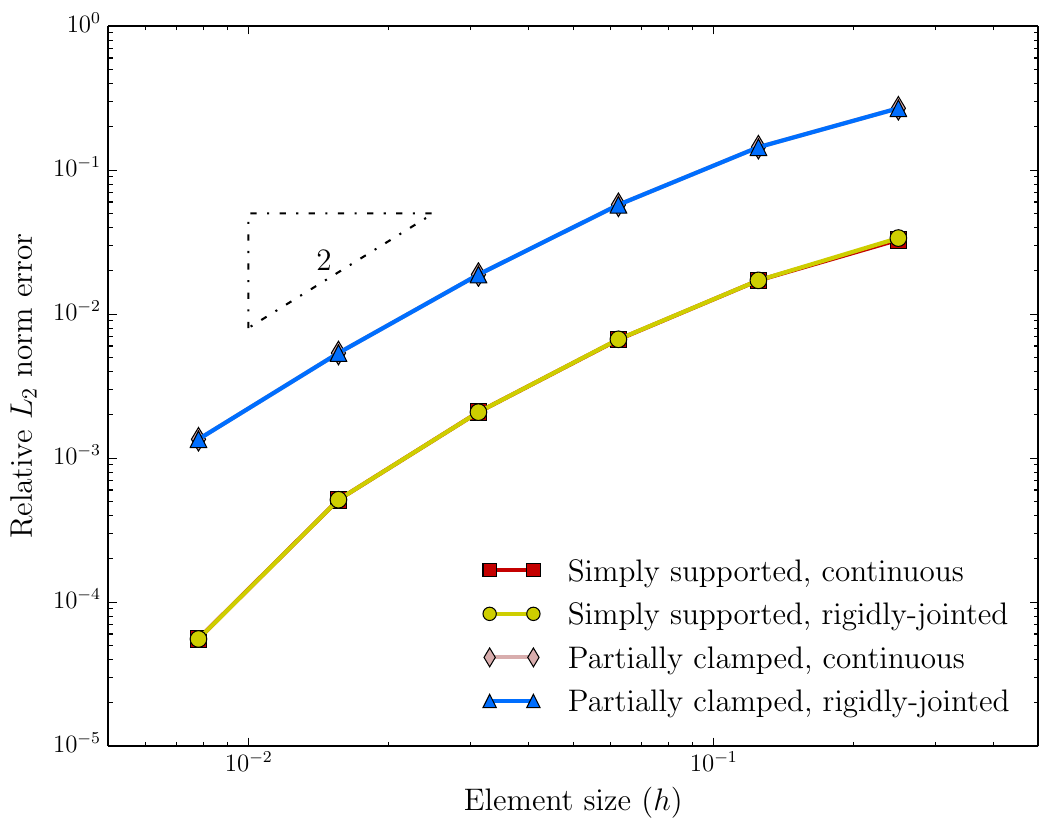}
  		\label{fig:plateErrorNonEV}
  	}
  	\subfloat[][Unstructured mesh] 
  	{
  		\includegraphics[scale=0.45]{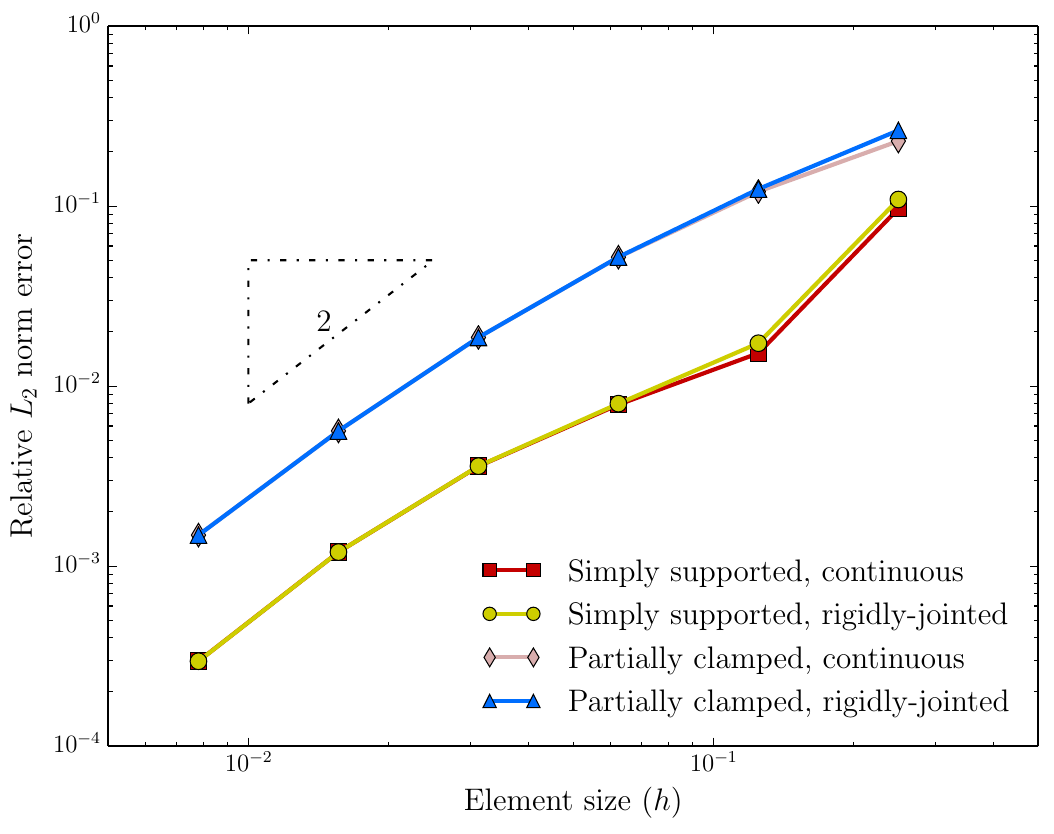}
  		\label{fig:plateErrorEV}
 	} 
	\caption{Convergence of the relative $L_2$ norm error for the simply supported and the partially clamped square plates. The errors for the continuous and rigidly-jointed plates are almost identical.}
	\label{fig:plateError}
\end{figure}

In Figure~\ref{fig:plateError} the convergence of the relative $L_2$ norm error with decreasing mesh size is plotted. In all cases the errors converge with the optimal rate of two.  Figure~\ref{fig:plateErrorNonEV} depicts the convergence for the structured mesh and Figure~\ref{fig:plateErrorEV}  for the unstructured mesh. Both plots contain the convergence of the fully simply supported plate as well as the partially clamped plate. In addition, the convergence of the plates with rigidly-jointed hinges are included. The continuous and rigidly-jointed plates yield the same results providing an evidence for the soundness of the proposed weak enforcement of rigid joint constraints. The errors for the less constrained fully simply supported plate are smaller than the errors for the partially clamped plate. It is noteworthy that the manifold-based basis functions can achieve the optimal convergence order even on an unstructured mesh. Moreover, the magnitude of the relative errors for the structured and unstructured meshes is almost the same.

%

\subsection{Pinched square tube \label{sec:tube}}
\begin{figure}[]
	\centering 
	\includegraphics[scale=1.0]{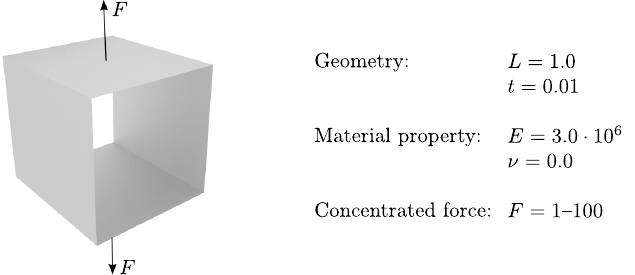}
	\caption{Geometry and loading of the pinched square tube. }
	\label{fig:tubeLoading}
\end{figure}

\begin{figure}[]
	\centering 
	\includegraphics[scale=0.275]{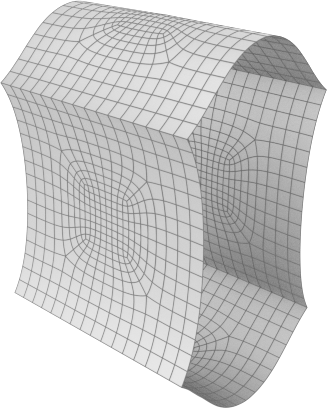}
	\hfil
	\includegraphics[scale=0.275]{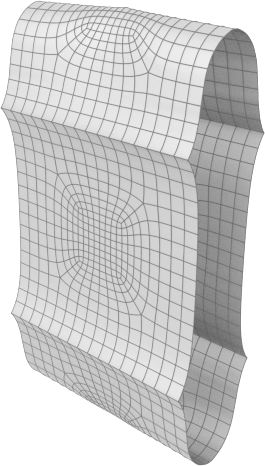}
	\hfil
	\includegraphics[scale=0.275]{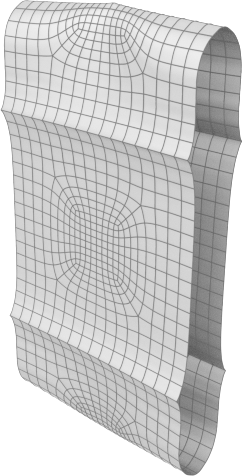}
	\caption{Deflected shapes of the pinched tube at load values $F=10$, $F= 50$ and $F=100$ (from left to right).}
	\label{fig:deformedTube}
\end{figure}

\begin{figure}[!ht]
	\centering 
	\includegraphics[scale=0.5]{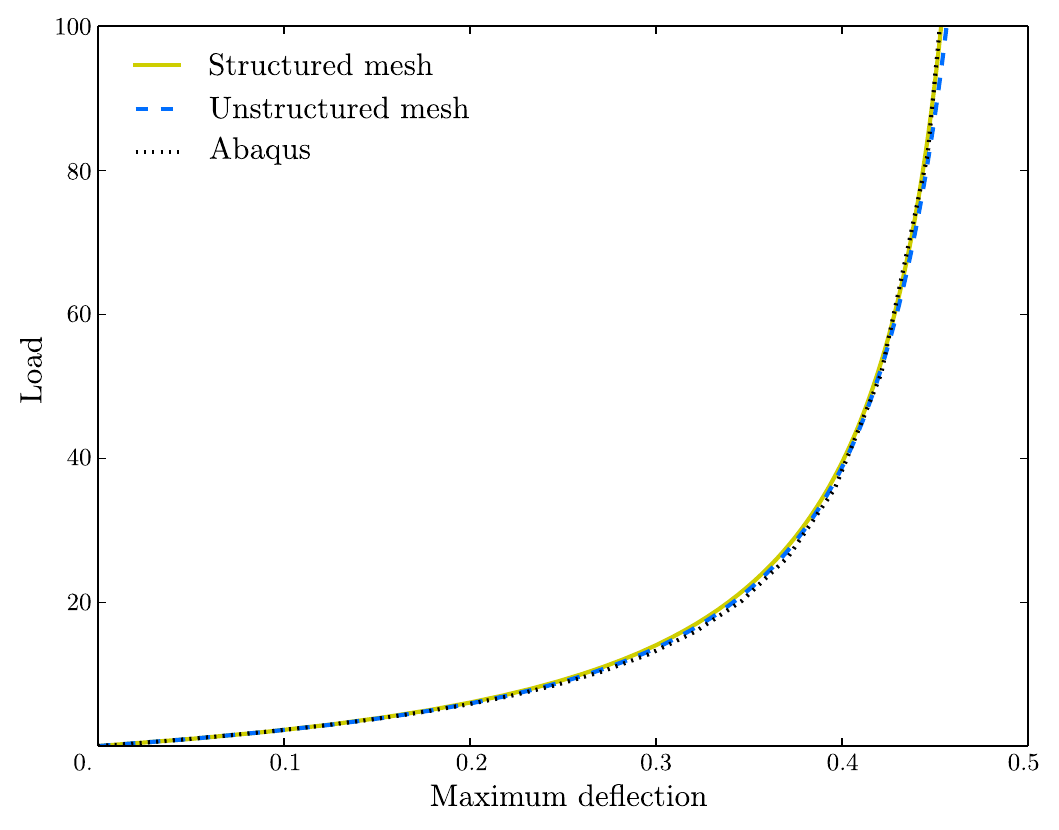}
	\caption{Load-displacement curve of the pinched square tube for the load attachment point.}
	\label{fig:pinchTube}
\end{figure}

%
This last example involves the geometrically nonlinear analysis of a pinched square tube subjected to two diametrically opposite concentrated forces, see Figure~\ref{fig:tubeLoading}. The tube consists out of two horizontal and two vertical plates that are rigidly connected along their edges. It is discretised either with a uniform structured mesh with an element size $h=L/16$ or with an unstructured mesh shown in Figure~\ref{fig:deformedTube}. The structured  and unstructured meshes have~4352 and~6400 nodes, respectively. To simulate the rigid joints between the plates the relevant edges in the mesh are tagged as a crease and the angle between the two normals across a crease are constrained to remain constant during deformation. This example is an adaptation of the pinched cylinder benchmark example widely used for comparing the performance of shell elements. The pinched square tube has also been previously considered in~\cite{Cirak:2011aa}. The deflected shapes at three different load values are shown in Figure~\ref{fig:deformedTube}.  As can be inferred from the deflected shapes, the tube exhibits large membrane deformations and localised bending deformations around the rigid creases and two load application points. Furthermore, it is evident that the right angle between the plates is preserved during deformation. In the load-displacement curve depicted in Figure~\ref{fig:pinchTube} the results for the structured and unstructured meshes are visually indistinguishable. Both are in close agreement with the reference solution obtained with the commercial finite element software Abaqus.  The Abaqus result is for a structured fine mesh with 10612 nodes and a six-parameter Reissner-Mindlin type shell theory.


 %

%% file: conclusions.tex
\section{Conclusions}                          
%
We have presented, by extending Majeed and Cirak~\cite{majeedCirak:2016},  new manifold-based basis functions for isogeometric design and analysis of surfaces with arbitrary smoothness, prescribed sharp features and boundaries.  The surface is described with a quadrilateral control mesh and~$C^0$ continuous creases are introduced along element edges tagged by the user. Manifold techniques are extremely versatile in the sense that the chart domains, their respective transition functions, local polynomial approximants, and blending functions can all be chosen to fit the needs of the specific application at hand. We introduce creases by choosing polynomial approximants that are piecewise $C^0$ continuous along the element edges and by modifying the geometry of the chart domains.  Due to the similarities between the $C^0$ continuous creases and boundaries, the new basis functions also simplify the treatment of boundaries. That is, different from~\cite{majeedCirak:2016}, the boundaries can be described without introducing an additional layer of ghost elements outside the domain. The chart domains have the shape of polygonal disks with curved boundaries and are composed out of quasi-conformally mapped unit squares representing the reference finite elements. The respective transition functions are easily computed by expressing the parametric coordinates of a unit square as a complex number. Furthermore, a new type of blending function is used which is assembled from tensor-product cubic B-spline pieces. In contrast to the ones in~\cite{majeedCirak:2016}, the new blending functions do not require normalisation so that the partition of unity approximation leads to non-rational basis functions. Finally, in order to obtain a mesh-based approximation scheme, the coefficients of the local polynomials are expressed as vertex coefficients using a least-squares procedure. The degree of the local polynomials has to be chosen such that there are no more coefficients than vertices in a chart domain. The number of vertices in a chart domain can be increased by introducing more vertices, i.e. refining the control mesh, while keeping the chart domain size constant. The obtained basis functions have a closed form analytic description, are locally supported and are polynomial in regular regions of the mesh. 

In closing, we note that the introduced manifold-based basis function construction may be interpreted as the extension of the conventional partition of unity method~\cite{melenk1996partition} to manifold surfaces. As detailed in~\cite{zhang2019manifold}, depending on the choice of the local polynomials and blending functions univariate manifold-based basis functions can be made to either reproduce B-splines or are identical to B-splines. Their convergence properties can be deduced from standard partition of unity results~\cite{melenk1996partition, babuvska1997partition}. Hence, manifold-based constructions may serve as a bridge  between isogeometric analysis using splines~\cite{Hughes:2005aa}  and the many partition of unity method inspired discretisation techniques, such as the generalised finite element method~\cite{strouboulis2000design}, hp-clouds~\cite{duarte1996h} or the extended finite element method~\cite{moes1999finite}; see also the recent review~\cite{chen2017meshfree}.  Specifically, manifold-based constructions can facilitate the consideration of industrial CAD geometries, in form of NURBS and other spline representations, in the partition of unity methods.  In turn, the very many enrichment techniques developed for partition of unity methods can be applied to isogeometric analysis by utilising them as local polynomial approximants in the manifold construction.  The exploration of these links suggests itself as a promising direction for future research.

%% file: appendix.tex
\appendix
\section{Conformal maps \label{sec:conformApp}} 
\begin{figure}[]
  \centering
  \subfloat[][Complex plane\label{fig:complexPlane1}] 
  {
	\includegraphics[scale=0.8]{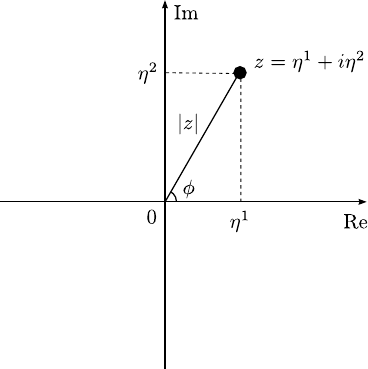}
  }\hspace{0.16\textwidth}
  \subfloat[][Rotation and scaling in the complex plane\label{fig:complexPlane2}] 
  {
	\includegraphics[scale=0.8]{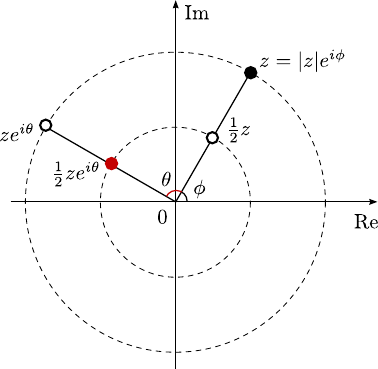}
  } 
  \caption{Illustration of a complex number in the complex plane and the operations of rotation and scaling.}
  \label{fig:complexplane}
\end{figure}

%
In this section we briefly motivate the quasi-conformal map used throughout this paper; see also~\cite{needham1998visual} for a general visually orientated introduction to complex analysis. In contrast to the conformal map used in Majeed and Cirak~\cite{majeedCirak:2016} the quasi-conformal map is not angle-preserving, but it is infinitely smooth except at the extraordinary vertex. With a conformal map an infinitesimal circle is mapped to a circle whereas with a quasi-conformal map it is mapped to an ellipse. The quasi-conformal map is essential for the introduced new chart domains, which are necessary for the construction of the creased basis functions. As discussed, three different quasi-conformal maps~$\vec{\varPsi}^{\textrm{I}}_j$,~$\vec{\varPsi}^{\textrm{II}}_j$ and~$\vec{\varPsi}^{\textrm{III}}_j$ are needed to map the reference element~$\Box$ onto elements in the chart domains~$\hat \Omega_j$ according to 
\begin{equation*}
	\left \{\vec{\varPsi}^{\textrm{I}}_j, \, \vec{\varPsi}^{\textrm{II}}_j, \, \vec{\varPsi}^{\textrm{III}}_j \right \} \colon  \vec \eta = (\eta^1, \, \eta^2) \in \Box \mapsto \vec \xi = (\xi_j^1, \, \xi_j^2) \in \hat \Omega_j \, .
\end{equation*}
In the complex plane the reference element coordinates~$\vec \eta = (\eta^1, \, \eta^2 )$ can be, as depicted in~Figure~\ref{fig:complexPlane1},  expressed either in Cartesian or polar form
\begin{equation*}
	z = \eta^1 + i\eta^2 = |z|(\cos \phi + i\sin \phi)  \, 
\end{equation*}
with the radius and the phase
\begin{equation*}
	|z| = \sqrt{(\eta^1)^2 + (\eta^2)^2}\,, \quad \phi = \arctan(\eta_2/\eta_1) \, .
\end{equation*}
With the Euler's formula the polar form can be equivalently expressed by 
\begin{equation*}
	z =  |z|e^{i\phi}  \, .
\end{equation*}
Considering a complex number as a vector from the origin to the point with the coordinates~$z$, a multiplication by a scalar represents a scaling and a multiplication by~$e^{i \theta}$ represents a rotation. For instance, Figure~\ref{fig:complexPlane2} illustrates that~$ze^{i\theta}$ is obtained by rotating~$z$ anticlockwise by an angle~$\theta$ and that~$1/2z$ is obtained by scaling the radius of~$z$ by one half. 

\begin{figure}[]
\centering
  \subfloat[][Conformal map with $\beta=4/v_j$ \label{fig:conformal1}] 
  {
	
	\includegraphics[scale=0.825]{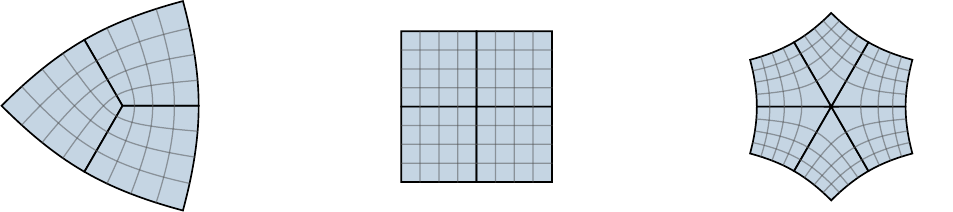}
  }
  \\
  \subfloat[][Quasi-conformal map with  $\beta=1$ \label{fig:conformal2}] 
  {
	\includegraphics[scale=0.825]{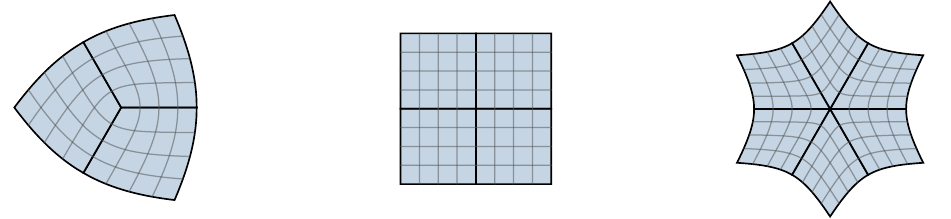}
  } 
  \caption{Comparison of the iso-parameter lines for the conformal  and the quasi-conformal map  depending on the value of the parameter $\beta$ for valences $v \in \{ 3,4,6\}$.
}
  \label{fig:conformal}
\end{figure}

In the quasi-conformal map introduced in~\eqref{eq:conformalMap}, repeated here for convenience, 
\begin{equation*}
	\vec{\Lambda^\textrm{I}_{qr}}(z;v_j,n_j) =  \frac{ |z|^{\beta}}{|z|^{4/v_j}}   z^{4/v_j}  e^{i 2\pi (n_j-1)/v_j} = |z|^\beta e^{i\phi 4/v_j}  e^{i 2\pi (n_j-1)/v_j} 	\, ,
\end{equation*}
the parameter $\beta$ controls how the radius of~$z$ is scaled.  For $\beta=4/v_j$  the mapping reduces to the conformal map used in~\cite{majeedCirak:2016}, i.e., 
\begin{equation*}
	\vec{\Lambda^\textrm{I}_{qr}}(z;v_j,n_j) =   z^{4/v_j}  e^{i 2\pi (n_j-1)/v_j} = |z|^{4/v_j}  e^{i 4 \phi /v_j }e^{i 2\pi (n_j-1)/v_j}  \, ,
\end{equation*}
which scales the radius of~$z$ to~$|z|^{4/v_j}\,$.  In the present paper we choose~$\beta=1$  so that the radius is not scaled. The choice of~$\beta$ on the conformal maps is illustrated in Figure~\ref{fig:conformal}. As visible, for~$\beta=1$ the parameter lines are not orthogonal to each other within the elements and the map is not angle-preserving.  Therefore, it is called a quasi-conformal map. For the introduced maps~$\vec{\varPsi}^{\textrm{II}}_j$ and~$\vec{\varPsi}^{\textrm{III}}_j$, i.e. ~\eqref{eq:conformalMap2} and~\eqref{eq:concavemap}, it is important that the radius of~$z$ is not scaled. This is necessary because each creased sector may have different number of elements. Not scaling the radius ensures that  the quasi-conformal map and its derivatives across crease edges are continuous c.f.~Figures \ref{fig:v5c2map}, \ref{fig:v5c3map} and \ref{fig:concavemap}. 
%

\section{Blending functions \label{sec:weightsApp}}
\begin{figure}[]
  \centering 
  \includegraphics[scale=0.875]{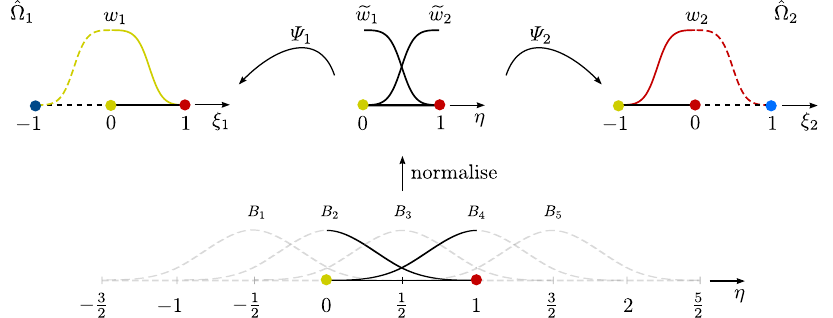}
  \caption{Construction of smooth rational blending functions from cubic B-splines defined over a parameter space with a knot-distance~$1/2$.}
  \label{fig:weightFuns}
\end{figure}

%
The blending functions proposed in this paper are polynomial in contrast to the rational blending functions used in Majeed and Cirak~\cite{majeedCirak:2016}.  The difference between the two constructions is best understood by comparing the two Figures~\ref{fig:1dWeightFun} and~\ref{fig:weightFuns}.  In Figure~\ref{fig:weightFuns} the construction process of the rational blending functions is illustrated. Notice that the cubic B-spline basis used within the reference element is not complete so that it does not add up to one. Therefore, normalisation is necessary to  satisfy the partition of unity property, resulting in rational blending functions. With the numbering introduced in Figure~\ref{fig:weightFuns} the rational blending functions are given by 
\begin{equation*}
	\widetilde{w}_j(\eta) = \frac{B_{2j}(\eta)}{\sum_{k=1}^{2} B_{2k}(\eta)} \quad \text{with} \quad w_j(\xi_j) = \widetilde w_j (\varPsi^{-1}_j (\xi_j))  \quad \text{and}  \quad j\in\{1,2\} \, .
\end{equation*}
 It is straightforward to extend this construction to the bivariate case. As a final remark, the rational blending functions have one knot and the polynomial blending functions used in this paper have three knots within the reference element. 
